\documentclass{article}

\usepackage{amsfonts,pstricks}
\usepackage{amsmath, amsthm, amssymb,  mathrsfs, epsfig}
\usepackage{hyperref}
\usepackage{graphicx}
\usepackage{tikz}
\usepackage{caption}
\usepackage{authblk}
\usepackage{bm}
\usepackage{tensor}
\usetikzlibrary{decorations.markings}
\usepackage{enumitem}
\hypersetup{linktocpage = true}

\DeclareMathOperator{\Hom}{Hom}
\DeclareMathOperator{\End}{End}

\DeclareMathOperator{\Tr}{Tr}

\def\Real{\mathbb{R}}

\def\Complex{\mathbb{C}}
\def\CatC{\mathcal{C}}
\def\CatD{\mathcal{D}}
\def\CatT{\mathcal{T}}

\def\CatG{\mathcal{G}}

\def\TQFT{\mathrm{TQFT}}

\def\unit{\mathbf{1}}
\def\CatG{\mathcal{G}}
\def\IC#1{\overline{#1}^{p}}
\def\DSC#1{\overline{#1}^{\oplus}}
\def\Sphere{\mathbb{S}}
\def\Label#1{L(#1)}
\def\Special{\text{special}}
\def\TVBW{\text{TVBW}}
\def\SET{\text{SET}}

\def\SMFC{\text{SMFC}}

\def\lsup#1#2{\tensor[^{#1}]{{#2}}{}}

\newtheorem{theorem}{Theorem}[section]

\newtheorem{example}[theorem]{Example}
\newtheorem{corollary}[theorem]{Corollary}
\newtheorem{proposition}[theorem]{Proposition}
\newtheorem{remark}[theorem]{Remark}
\newtheorem{definition}[theorem]{Definition}

\def \Yshape #1#2#3#4#5
{
\draw (#1,-1+ #2)--(#1,#2) -- (-1+ #1,1 + #2);
\draw (#1,#2) -- (1+ #1,1+ #2);
\draw (-1.5+#1,1+ #2) node{{\tiny ${#3}$}};
\draw (1.5+#1, 1+ #2) node{{\tiny ${#4}$}};
\draw (-0.5+#1,-1+ #2) node{{\tiny ${#5}$}};
}
\def \InverseYshape #1#2#3#4#5
{
\draw (#1,1+ #2)--(#1,#2) -- (-1+ #1,-1 + #2);
\draw (#1,#2) -- (1+ #1,-1+ #2);
\draw (-1.5+#1,-1+ #2) node{{\tiny ${#3}$}};
\draw (1.5+#1, -1+ #2) node{{\tiny ${#4}$}};
\draw (-0.5+#1,1+ #2) node{{\tiny ${#5}$}};
}

\def \Phishape #1#2#3#4#5#6
{
\Yshape{#1}{#2}{#3}{#4}{#5}
\begin{scope}[yshift = 2cm]
\draw (#1,1+ #2)--(#1,#2) -- (-1+ #1,-1 + #2);
\draw (#1,#2) -- (1+ #1,-1+ #2);
\draw (-0.5+#1,1+ #2) node{{\tiny ${#6}$}};
\end{scope}
}

\def \Ishape #1#2#3#4#5#6#7
{
\Yshape{#1}{2+#2}{#5}{#6}{#7}
\draw (#1,1+ #2)--(#1,#2) -- (-1+ #1,-1 + #2);
\draw (#1,#2) -- (1+ #1,-1+ #2);
\draw (-1.5+#1,-1+ #2) node{{\tiny ${#3}$}};
\draw (1.5+#1, -1+ #2) node{{\tiny ${#4}$}};
}

\def \ShortIshape #1#2#3#4#5#6#7
{
\draw (#1,1+ #2)--(#1,#2) -- (-1+ #1,-1 + #2);
\draw (#1,#2) -- (1+ #1,-1+ #2);
\draw (-1+#1, 2+#2) -- (#1, 1+#2) -- (1+#1, 2+#2);
\draw (-1.5+#1,-1+ #2) node{{\tiny ${#3}$}};
\draw (1.5+#1, -1+ #2) node{{\tiny ${#4}$}};
\draw (-1.5+#1, 2+ #2) node{{\tiny ${#5}$}};
\draw (1.5+#1, 2+ #2) node{{\tiny ${#6}$}};
\draw (-0.5+#1, 0.5+ #2) node{{\tiny ${#7}$}};
}

\def \IshapeL #1#2#3#4#5
{
\draw (#1,1+ #2)--(#1,#2) -- (-1+ #1,-1 + #2);
\draw (#1,#2) -- (1+ #1,-1+ #2);
\draw (-1+#1, 2+#2) -- (#1, 1+#2) -- (1+#1, 2+#2);
\draw (-0.8+#1, 0.5+ #2) node{{\tiny ${#3}$}};
\draw (#1,-0.3+#2) node{{\tiny $#4$}};
\draw (#1,1.3+#2) node{{\tiny $#5$}};
}

\def \IshapeR #1#2#3#4#5
{
\draw (#1,1+ #2)--(#1,#2) -- (-1+ #1,-1 + #2);
\draw (#1,#2) -- (1+ #1,-1+ #2);
\draw (-1+#1, 2+#2) -- (#1, 1+#2) -- (1+#1, 2+#2);
\draw (0.8+#1, 0.5+ #2) node{{\tiny ${#3}$}};
\draw (#1,-0.3+#2) node{{\tiny $#4$}};
\draw (#1,1.3+#2) node{{\tiny $#5$}};
}

\def \IshapeRP #1#2#3#4#5#6
{
\draw (#1,1+ #2)--(#1,#2) -- (-1+ #1,-1 + #2);
\draw (#1,#2) -- (1+ #1,-1+ #2);
\draw (-1+#1, 2+#2) -- (#1, 1+#2) -- (1+#1, 2+#2);
\draw (0.6+#1, 0.5+ #2) node{{\tiny ${#3}$}};
\draw [dotted] plot [smooth] coordinates {(-0.6+#1,-0.3+#2) (-0.9+#1,0.5+#2) (-0.6+#1,1.3+#2)};
\draw [dotted] plot [smooth] coordinates {(0.9+#1,-0.3+#2) (1.2+#1,0.5+#2) (0.9+#1,1.3+#2)};
\draw (-1.1+#1,1.1+#2) node{{\tiny $#6$}};
\draw (#1,-0.3+#2) node{{\tiny $#4$}};
\draw (#1,1.3+#2) node{{\tiny $#5$}};
}

\def \Vline #1#2#3#4
{
\draw (#1,#2) -- (#1,#2 + #3);
\draw (-0.5+ #1,#2) node{\tiny ${#4}$};
}

\def \VlineDot #1#2#3
{
\draw (#1, #2) -- (#1, #2 + #3);
}

\def \Braid#1#2
{
\draw (2+#1,#2) -- (#1,2+#2);
\draw (0.8+#1,0.8+#2)[color = white, line width = 2mm] -- (1.2+#1,1.2+#2);
\draw (#1,#2) -- (2+#1,2+#2);
}
\def \InvBraid#1#2
{
\draw (#1,#2) -- (2+#1,2+#2);
\draw (2+#1,#2)[color = white, line width = 2mm] -- (#1,2+#2);
\draw (2+#1,#2) -- (#1,2+#2);
}

\def\RecMor#1#2#3#4#5#6 
{
\draw[postaction={decorate}] (#1,0.5+#2) -- (#1,#2) node[left]{{\tiny $#3$}};
\draw[postaction={decorate}] (1+#1,0.5+#2) -- (1+#1,#2)node[right]{{\tiny $#4$}};

\draw[postaction={decorate}] (#1,2+#2)node[left]{{\tiny $#5$}} -- (#1,1.5+#2);
\draw[postaction={decorate}] (1+#1,2+#2)node[right]{{\tiny $#6$}} -- (1+#1,1.5+#2);

\draw (-0.5+#1,0.5+#2) rectangle (1.5+#1,1.5+#2);
}

\def \Segment#1#2#3#4#5
{
\tikz[scale = 0.5, baseline = {(current bounding box.center)}]{
 \draw (#1,#2) -- (2+ #1,#2);
 \draw (-0.5+#1,#2) node{{\tiny $#3$}};
 \draw (2.5+#1,#2) node{{\tiny $#4$}};
 \draw (1+#1,0.5+#2) node{{\tiny $#5$}};
}
}
\def \Triangle#1#2#3#4#5#6#7#8#9
{
\tikz[scale = 0.5, baseline= {(current bounding box.center)}]{
 \draw (#1,#2) -- (2+#1,#2) -- (1+#1,1.73+#2) -- (#1,#2);
 \draw (-0.5+#1,#2) node{{\tiny $#3$}};
 \draw (2.5+#1,#2) node{{\tiny $#4$}};
 \draw (1+#1,2.23+#2) node{{\tiny $#5$}};
 \draw (1+#1,-0.3+#2) node{{\tiny $#6$}};
 \draw (2+#1,1+#2) node{{\tiny $#7$}};
 \draw (0+#1,1+#2) node{{\tiny $#8$}};
 \draw (1+#1,0.5+#2) node{{\tiny $#9$}};
}
}

\begin{document}

\title{State Sum Invariants of Three Manifolds from Spherical Multi-fusion Categories}
\author{Shawn X. Cui%
       \thanks{Email: \texttt{xingshan@stanford.edu}}}
\affil{Stanford Institute for Theoretical Physics,\\ Stanford University, Stanford, California 94305}
\author{Zhenghan Wang%
       \thanks{Email: \texttt{zhenghwa@microsoft.com}}}
\affil{Microsoft Station Q and Dept. of Math., \\University of California, Santa Barbara, California 93106}
\maketitle

\begin{abstract}
We define a family of quantum invariants of closed oriented $3$-manifolds using spherical multi-fusion categories. The state sum nature of this invariant leads directly to $(2+1)$-dimensional topological quantum field theories ($\TQFT$s), which generalize the Turaev-Viro-Barrett-Westbury ($\TVBW$) $\TQFT$s from spherical fusion categories. The invariant is given as a state sum over labeled triangulations, which is mostly parallel to, but richer than the $\TVBW$ approach in that here the labels live not only on $1$-simplices but also on $0$-simplices. It is shown that a multi-fusion category in general cannot be a spherical fusion category in the usual sense. Thus we introduce the concept of a spherical multi-fusion category by imposing a weakened version of sphericity. Besides containing the $\TVBW$ theory, our construction also includes the recent higher gauge theory $(2+1)$-$\TQFT$s  given by Kapustin and Thorngren,  which was not known to have a categorical origin before.
\end{abstract}

\section{Introduction}
\label{sec:intro}
A fundamental connection between three dimensional topology and higher categories is the $(2+1)$-dimensional topological quantum field theory ($\TQFT$) introduced in \cite{witten1988topological,atiyah1988topological}. A $(2+1)$-$\TQFT$ associates to every
closed oriented $2$-manifold a finite dimensional vector space and to every $3$-manifold a vector in the vector
space corresponding to its boundary. These assignments should satisfy certain axioms \cite{atiyah1988topological}. The empty set is considered as a closed $2$-manifold and the vector space
associated to it is required to be $\Complex$. Then the vector corresponding to a closed $3$-manifold becomes a complex scalar called the partition function or path integral, which is an invariant of $3$-manifolds. Invariants arising from $\TQFT$s are called quantum invariants.

Quantum invariants have largely been constructed by state-sum models from monoidal categories and Hopf algebras. Reshetikhin and Turaev constructed an invariant of $3$-manifolds using modular tensor categories, which is believed to be the mathematical realization of Witten's $\TQFT$ from non-abelian Chern-Simon theories \cite{reshetikhin1991invariants}. Turaev and Viro gave a state-sum invariant of $3$-manifolds (Turaev-Viro invariant) from a ribbon fusion category \cite{turaev1992state}. Later Barrett and Westbury generalized this construction (Turaev-Viro-Barrett-Westbury invariant or $\TVBW$ invariant) by using spherical fusion categories \cite{barrett1996invariants}. These invariants can all be extended to define a $(2+1)$-$\TQFT$. Apart from these categorical constructions, another approach is by using certain Hopf algebras, among which the Kuperberg invariant \cite{kuperberg1996noninvolutory} and the Hennings invariant \cite{kauffman1995invariants} \cite{hennings1996invariants} are non-semisimple generalizations of the Turaev-Viro invariant and the Reshetikhin-Turaev invariant, respectively. A special case of the Kuperberg invariant (and also the Turaev-Viro invariant) reduces to the Dijkgraaf-Witten theory \cite{dijkgraaf1990topological}. The study of $(2+1)$-$\TQFT$s has led to applications in quantum groups, $3d$ topology, and knot theories. For example, the Turaev-Viro invariant can distinguish certain $3$-manifolds which are homotopy equivalent.

The main result of this paper is a construction of a state-sum invariant of $3$-manifolds from what we call a spherical multi-fusion category ($\SMFC$). When the $\SMFC$ is a fusion category, the invariant reduces to the $\TVBW$ invariant. It is straightforward to extend the construction to obtain a $(2+1)$-$\TQFT$. However, for simplicity, here we only focus on quantum invariants. Our contribution touches on the following three aspects.

Firstly, we introduced the concept of $\SMFC$s. The current definitions of spherical categories and multi-fusion categories are in general not compatible; a multi-fusion category can never be spherical unless it is a fusion category (see Section \ref{subsec:multifusion}). Thus a $\SMFC$ is not a spherical category in the usual sense. We weakened the definition of sphericity based on the construction of state-sum invariants. Explicitly, let $\CatC = \bigoplus\limits_{i,j \in I} \CatC_{ij}$ be a multi-fusion category, where $\CatC_{ij}$, called the $(i,j)$-sector, satisfies $\CatC_{ik} \otimes \CatC_{kj} \subset \CatC_{ij}$. Here $I$ is called the index set and for each $i \in I$, $\CatC_{ii}$ is a fusion category with unit $\unit_i$ and the unit of $\CatC$ is $\unit = \bigoplus\limits_{i\in I}\unit_i$. If $f \in \End(X), X \in \CatC_{ij}$, instead of requiring the left trace of $f$ equal the right trace of $f$ on the nose, i.e., $\Tr^l(f) = \Tr^r(f)$ which in general does not hold, we define $\CatC$ to be spherical if $|\Tr^l(f)| = |\Tr^r(f)|$ where $|\Tr^l(f)|$, $|\Tr^l(f)|$ are scalars that satisfy $\Tr^l(f) = |\Tr^l(f)| id_{\unit_j}$, $\Tr^r(f) = |\Tr^r(f)| id_{\unit_i}$. When $I$ consists of one element, this definition reduces to the usual one. Another motivation of this weakening comes from graphical calculus of multi-fusion categories. The new definition of sphericity guarantees that isotopic colored graphs in the sphere have the same evaluation. 

Secondly, the construction of quantum invariants is more general than the $\TVBW$ approach. In the $\TVBW$ model, only the $1$-simplices are colored while in our model, both the $1$-simplices and $0$-simplices are colored. Let $M$ be a closed oriented $3$-manifolds and $\CatT$ be a triangulation of $M$ whose vertices are ordered. Let $\CatC$ be a $\SMFC$ with index set $I$. A coloring $F$ of $\CatT$ assigns to each $0$-simplex ordered by $i$ an element $f^0_i \in I$ and to each $1$-simplex $(ij)$ a simple object $f^1_{ij}\in \CatC_{f^0_if^0_j}$. Then the partition function is defined to be,
\begin{align}
\label{equ:partition2}
Z_{\CatC}(M,\CatT) = \sum\limits_{F = (f^0,f^1)} \prod\limits_{\tau \in \CatT^0} K^{-1} \prod\limits_{\alpha \in \CatT^1} d_{f^1_{\alpha}} \,  \prod\limits_{\beta \in \CatT^2} \theta_F(\beta)^{-1} \, \prod\limits_{\gamma \in \CatT^3} \tilde{Z}^{\epsilon(\gamma)}_F(\gamma;B),
\end{align}
where $K$, $d_{(\cdot)}$, $\theta_{F}(\cdot)$, and $\tilde{Z}_{F}^{\epsilon(\cdot)}(\cdot, B)$ are certain scalars associated to simplices of various dimensions. See Section \ref{sec:main} for details. The main result is as follows.

\begin{theorem}[Main]
\label{thm:main}
The partition function $Z_{\CatC}(M,\CatT)$ is independent of the triangulation $\CatT$, and thus $Z_{\CatC}(M):= Z_{\CatC}(M,\CatT)$ is an invariant of closed oriented $3$-manifolds. Moreover, this construction extends to a $(2+1)$-$\TQFT$.
\end{theorem}

Lastly, by studying a class of $\SMFC$s coming from categorical groups and some additional cohomological data, we recovered the $(2+1)$-$\TQFT$ in \cite{kapustin2013higher} which is obtained from higher gauge theory. The latter $\TQFT$ was not known to have a categorical construction.

The rest of the paper is organized as follows. In Section \ref{sec:review} we provide a review of basic category theories and propose the concept of $\SMFC$s. Section \ref{sec:main} contains the main construction of quantum invariants. In Section \ref{sec:generalized}, we define generalized categorical groups and study the $\SMFC$s obtained from them. Finally in Section \ref{sec:SET}, we make some connections to symmetry enriched topological phases.

\section{Spherical Multi-fusion Categories}
\label{sec:review}
We assume the readers to have a background on basic category theories and especially monoidal (tensor) category theories. In Section \ref{subsec:pivotal} we set up some notations and briefly review monoidal categories with additional structures such as duals and pivotal structures. We also review graphical calculus which is a convenient way to represent morphisms. There are many references on this subject. For instance, see \cite{bojko2001lectures}\cite{Etingof_finitetensor}\cite{turaev1994quantum}\cite{kassel1995quantum}, etc. In Section \ref{subsec:multifusion}, we first recall the definition of (multi)-fusion categories and then introduce the concept of a spherical multi-fusion category, which is not spherical according to the usual definition unless the category is a fusion category. Spherical multi-fusion categories are natural generalizations of spherical fusion categories.

Throughout the context, let $\CatC$ be a category. Denote by $\CatC^0$ the set of objects and by $\Hom_{\CatC}(X, Y)$ or simply $\Hom(X, Y)$ the set of morphisms between an object $X$ and an object $Y$. If $X=Y$, $\Hom(X,X)$ is also written as $\End(X)$. The compositions of morphisms will be read from right to left. Namely, if $f \in \Hom(X, Y), g \in \Hom(Y, Z)$, then $g\circ f \in \Hom(X, Z)$.

\subsection{Pivotal Categories and Graphical Calculus}
\label{subsec:pivotal}

Let $\CatC$ be a rigid monoidal category, that is, a category endowed with the tuple $(\otimes, \unit, a, l, r, (\cdot)^*)$, where $\otimes: \CatC \times \CatC \longrightarrow \CatC$ is the tensor product functor, $\unit$ is the unit object, and $a, l, r$ are natural isomorphisms:
\begin{align*}
a_{X,Y,Z}: & (X \otimes Y) \otimes Z \overset{\simeq}{\longrightarrow} X \otimes (Y \otimes Z), \\
l_X: &\unit \otimes X \overset{\simeq}{\longrightarrow}  X,\\
r_X: & X \otimes \unit \overset{\simeq}{\longrightarrow} X, \quad X, Y, Z \in \CatC^0,
\end{align*}
which satisfy the Pentagon Equation and Triangle Equation, and $(\cdot)^*$ is the contra-variant functor of taking duals. For each object $X$, denote the birth (also called co-evaluation) and death (also called evaluation) morphism by $b_X$ and $d_X$ respectively:
\begin{align*}
b_X: \unit \longrightarrow X \otimes X^*, \qquad d_X: X^* \otimes X \longrightarrow \unit,
\end{align*}


A pivotal structure on a rigid monoidal category is a natural isomorphism $\delta: Id_{\CatC} \longrightarrow (\cdot)^{**}$. Thus for each object $X$, there is an isomorphism:
\begin{align*}
\delta_{X}: X \overset{\simeq}{\longrightarrow} X^{**},
\end{align*}
such that $\delta_{X} \otimes \delta_{Y} \overset{\cdot}{=} \delta_{X \otimes Y}$, where \lq$\overset{\cdot}{=}$' means \lq equal' up to a composition of certain canonical isomorphism. With the pivotal structure $\delta$, we can define another set of \lq birth' and \lq death' morphisms,
\begin{align*}
b_{X}': \unit \longrightarrow X^* \otimes X, \qquad d_{X}': X \otimes X^* \longrightarrow \unit
\end{align*}
where $b_{X'}:= (id_{X^*} \otimes \delta_{X}^{-1})b_{X^*}$, $d_{X'}:= d_{X^*}(\delta_X \otimes id_{X^*})$.

A pivotal category is a rigid monoidal category with a chosen pivotal structure. Given a morphism $f \in \End(X)$, define the left trace $\Tr^l(f) \in \End(\unit)$ by $\Tr^l(f) = d_X(id_{X^*} \otimes f)b_X'$ and the right trace $\Tr^r(f) \in \End(\unit)$ by $\Tr^r(f) = d_{X}'(f \otimes id_{X^*})b_X$. Then a pivotal category is called spherical if $\Tr^l(f) = \Tr^r(f)$ for all endmorphisms $f$. Every pivotal category $\CatC$ is equivalent (in a properly defined sense) to a strict pivotal category $\hat{\CatC}$ where all the structural isomorphisms $a, l ,r, \delta $ are the identity map \cite{barrett1999spherical}\cite{ng2007higher}. If $\CatC$ is spherical, so is $\hat{\CatC}$.

In a strict pivotal category $\CatC$, graphical calculus is a convenient way to represent and manipulate morphisms. We sketch the rules for graphical calculus following the conventions in \cite{turaev1994quantum}\cite{turaev2010on}.

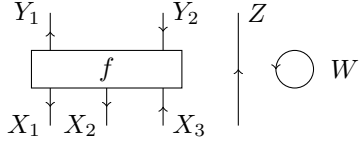
\begin{figure}
 \centering
\begin{tikzpicture}[scale = 0.5]
 \begin{scope}[decoration={
    markings,
    mark=at position 0.5 with {\arrow{>}}}]
    \draw [postaction={decorate}] (0,1) -- (0,0) node[left]{$X_1$};
    \draw [postaction={decorate}] (1.5,1) -- (1.5,0) node[left]{$X_{2}$};
    \draw [postaction={decorate}] (3,0)node[right]{$X_3$} -- (3,1) ;
    \draw [postaction={decorate}] (0,2) -- (0,3) node[left]{$Y_1$};
    \draw [postaction={decorate}] (3,3)node[right]{$Y_2$} -- (3,2) ;
    \draw (-0.5,1) rectangle (3.5,2)node[pos = 0.5]{$f$};

    \draw [postaction={decorate}] (5, 0) -- (5,3)node[right]{$Z$};
    \draw [postaction={decorate}] (6.5,1.5) circle (0.5cm) ;
    \draw (7.2, 1.5) node[right]{$W$};

    \end{scope}
 \end{tikzpicture}
 \caption{Example of a colored graph diagram. $f$ is a morphism from $X_1 \otimes X_2 \otimes X_3^*$ to $ Y_1^* \otimes Y_2$, and the whole diagram is a morphism from $X_1 \otimes X_2 \otimes X_3^* \otimes Z^*$ to $Y_1^* \otimes Y_2 \otimes Z^*$.}\label{fig:examplegraph}
\end{figure}

A graph diagram is a collection of rectangles \footnote{In \cite{turaev1994quantum}, they are called coupons.} and directed arcs (including circles) in $\Real \times [0,1]$, satisfying the conditions:
\begin{itemize}
\item The longer sides of each rectangle are parallel to $\Real \times \{0\}$ and the shorter sides vertical to $\Real \times \{0\}$.
\item The collection of $\{$rectangles, arcs$\}$ are mutually disjoint from each other except that every non-circular arc starts and ends transversely either on $\Real \times \{0,1\}$ or on the horizontal sides of a rectangle.
\end{itemize}
A $\CatC$-colored (or colored, for short) graph diagram is a graph diagram $\CatG$ which further satisfies:
\begin{itemize}
\item Each arc is labeled by an object of $\CatC$ and each rectangle is labeled by a morphism obeying the following rule. For a rectangle $\gamma$ labeled by $f$, denote by $\beta_1, \cdots, \beta_{m}$ the set of arcs incident to the bottom of $\gamma$, and by $\beta^1, \cdots, \beta^n$ the set of arcs incident to the top of $\gamma$ both listed from left to right. Note that some $\beta_{i_1}$ and $\beta_{i_2}$ might be the same arc if this arc intersects $\gamma$ twice. For each $\beta_i$ (resp. $\beta^j$), define $\epsilon_i$ (resp. $\epsilon^j$) to be $+1$ if $\beta_i$ (resp. $\beta^j$) is directed downwards near the rectangle, and $-1$ otherwise. Denote the label on the arc $\beta_i$ (resp. $\beta^j$) by $X_i$ (resp. $X^j$), then we require
\begin{align*}
f \in \Hom(\bigotimes\limits_{i=1}^{m} X_i^{\epsilon_i}, \bigotimes\limits_{j=1}^{n} (X^j)^{\epsilon^j}),
\end{align*}
where for an object $X$, $X^{+1}:= X$ and $X^{-1}:= X^*$. If $m=0$ or $n=0$, then let the corresponding object be the unit $\unit$. For instance, in Figure \ref{fig:examplegraph}, we represent the labels by putting an object next to each arc and a morphism inside each rectangle. Then $f \in \Hom(X_1 \otimes X_2 \otimes X_3^*,\, Y_1^* \otimes Y_2)$. Note that we have omitted and will never draw the lines $\Real \times \{0,1\}$.
\end{itemize}

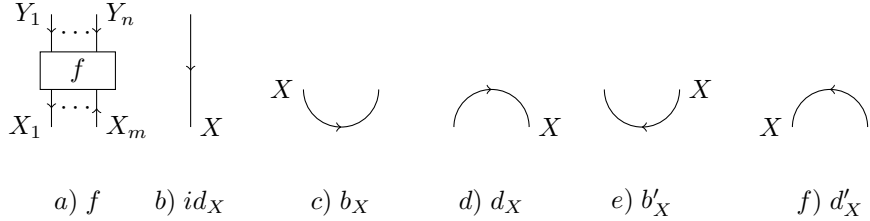
\begin{figure}
\centering
\begin{tikzpicture}[scale = 0.5]
 \begin{scope}[xshift = -0.5cm, decoration={
    markings,
    mark=at position 0.5 with {\arrow{>}}}]
 \draw [postaction={decorate}] (-0.2,1) -- (-0.2,0)node[left]{$X_1$};
 \draw [postaction={decorate}] (1,0)node[right]{$X_m$} -- (1,1);
 \draw (0.5, 0.5) node{$\cdots$};
 \draw (-0.5,1) rectangle (1.5, 2) node[pos = 0.5]{$f$};
 \draw [postaction={decorate}]  (-0.2,3) node[left]{$Y_1$} -- (-0.2,2);
 \draw [postaction={decorate}] (1,3)node[right]{$Y_n$} -- (1,2);
 \draw (0.5, 2.5) node{$\cdots$};
 \draw (0.5,-2) node{$a)\; f$};
 \end{scope}

 \begin{scope}[xshift = 3cm, decoration={
    markings,
    mark=at position 0.5 with {\arrow{>}}}]
 \draw [postaction={decorate}] (0,3) -- (0,0) node[right]{$X$};
 \draw (0,-2) node{$b)\; id_X$};
 \end{scope}

 \begin{scope}[xshift = 6cm, decoration={
    markings,
    mark=at position 0.5 with {\arrow{>}}}]
 \draw [postaction={decorate}] (0,1) node[left]{$X$} arc(-180:0:1cm) ;
 \draw (1,-2) node{$c)\; b_X$};
 \end{scope}

  \begin{scope}[xshift = 10cm, decoration={
    markings,
    mark=at position 0.5 with {\arrow{>}}}]
 \draw [postaction={decorate}]   (0,0)  arc(180:0:1cm) node[right]{$X$};
 \draw (1,-2) node{$d)\; d_X$};
 \end{scope}

 \begin{scope}[xshift = 14cm, decoration={
    markings,
    mark=at position 0.5 with {\arrow{>}}}]
 \draw [postaction={decorate}]   (2,1) node[right]{$X$} arc(0:-180:1cm);
 \draw (1,-2) node{$e)\; b_X'$};
 \end{scope}

  \begin{scope}[xshift = 19cm, decoration={
    markings,
    mark=at position 0.5 with {\arrow{>}}}]
 \draw [postaction={decorate}]   (2,0)  arc(0:180:1cm) node[left]{$X$};
 \draw (1,-2) node{$f)\; d_X'$};
 \end{scope}
\end{tikzpicture}
\caption{$F(\CatG)$ for some colored graph diagrams} \label{fig:rigid pic}
\end{figure}

Given a $\CatC$-colored graph diagram $\CatG$, denote by $\beta_1, \cdots, \beta_{m}$ (resp. $\beta^1, \cdots, \beta^{n}$) the set of arcs incident to $\Real \times \{0\}$ (resp. $\Real \times \{1\}$), and define the $X_i\,$'s, $\epsilon_i\,$'s, $X^j\,$'s, and $\epsilon^j\,$'s in the same way as above. Let $S(\CatG) = \bigotimes\limits_{i=1}^{m} X_i^{\epsilon_i}$ and $ T(\CatG) = \bigotimes\limits_{j=1}^{n} (X^j)^{\epsilon^j}$. Then $\CatG$ can be interpreted as a morphism $F(\CatG) \in \Hom(S(\CatG), T(\CatG))$ by the following rules:
\begin{enumerate}
\item If there is a color-preserving isotopy between $\CatG$ and $\CatG'$ relative to $\Real \times \{0,1\}$, then $F(\CatG) = F(\CatG')$.
\item If $\CatG$ is cut into two colored graph diagrams $\CatG_1$ (lower) and $\CatG_2$ (upper) by the line $\Real \times \{\frac{1}{2}\}$, then $F(\CatG) = F(\CatG_2)F(\CatG_1)$.
\item If $\CatG$ is separated into two disjoint colored graph diagrams $\CatG_1$ (left) and $\CatG_2$ (right) by the line $\{0\} \times \Real$, then $F(\CatG) = F(\CatG_1) \otimes F(\CatG_2)$.
\item Reversing the direction of an arc and changing its color to the dual at the same time does not change $F(\CatG)$.
\item If $\CatG$ is one of the diagrams in Figure \ref{fig:rigid pic}, then $F(\CatG)$ is the corresponding morphism listed below.
\end{enumerate}

It is not hard to see the above rules uniquely determine $F(\CatG)$. However, it takes more effort to check these rules are consistent. See \cite{turaev2010on} for more details.


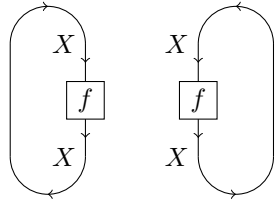
\begin{figure}
\centering
\begin{tikzpicture}[scale = 0.5]
 \begin{scope}[decoration={
    markings,
    mark=at position 0.5 with {\arrow{>}}}]
    \draw [postaction={decorate}] (2,1) arc(0:-180:1cm);
    \draw [postaction={decorate}] (0,4) arc(180:0:1cm);
    \begin{scope}[xshift = 2cm, yshift = 1cm]
      \draw [postaction={decorate}] (0,1) -- (0,0)node[left]{$X$};
      \draw (-0.5,1) rectangle (0.5, 2) node[pos = 0.5]{$f$};
      \draw [postaction={decorate}]  (0,3) node[left]{$X$} -- (0,2);
    \end{scope}
    \draw (0,1) -- (0,4);
 \end{scope}

 \begin{scope}[xshift = 5cm, decoration={
    markings,
    mark=at position 0.5 with {\arrow{>}}}]
    \draw [postaction={decorate}] (0,1) arc(-180:0:1cm);
    \draw [postaction={decorate}] (2,4) arc(0:180:1cm);
    \begin{scope}[xshift = 0cm, yshift = 1cm]
      \draw [postaction={decorate}] (0,1) -- (0,0)node[left]{$X$};
      \draw (-0.5,1) rectangle (0.5, 2) node[pos = 0.5]{$f$};
      \draw [postaction={decorate}]  (0,3) node[left]{$X$} -- (0,2);
    \end{scope}
    \draw (2,1) -- (2,4);
 \end{scope}
 \end{tikzpicture}
 \caption{$\Tr^l(f)$ (Left) and $\Tr^r(f)$ (Right)}\label{fig:trace}
\end{figure}

 A graph diagram is called closed if it does not have any free ends, i.e., it is disjoint from $\Real \times \{0,1\}$. If $\CatG$ is closed, then $F(\CatG) \in \Hom(\unit,\unit)$. We can view closed colored graph diagrams as sitting in $\Real^2$, and isotopy in $\Real^2$ does not change its value. For instance, given $f \in \End(X)$, the left and right trace of $f$ are represented by the two closed colored graph diagrams in Figure \ref{fig:trace}. Note that these two diagrams are not isotopic in $\Real^2$. However they become isotopic when we embed $\Real^2$ in $\Sphere^2 = \Real^2 \cup \{\infty\}$. Thus, a necessary condition for isotopic closed graph diagrams in $\Sphere^2$ to represent the same morphism is $\Tr^l = \Tr^r$, i.e., $\CatC$ is spherical. In fact, this condition is also sufficient. Let $\CatG$ be a closed colored graph diagram in $\Sphere^2$ (defined similarly as above), then we can remove a point $pt$ in the complement of $\CatG$, and view $\CatG$ as in $\Real^2 = \Sphere^2 \setminus \{pt\}$, and interpret $\CatG$ as a morphism $F(\CatG) \in \End(\unit)$. In \cite{turaev2010on} it is shown that if $\CatC$ is spherical, then $F(\CatG)$ is well-defined for $\CatG$ in $\Sphere^2$. 

If $\CatC$ is spherical, then define $\Tr(f):= \Tr^l(f) = \Tr^r(f).$ In particular, for an object $X$, define the dimension of $X$ to be $d_X:= \Tr(id_X)$, which is represented by a circle labeled by $X$. The direction of the circle is irrelevant since $\CatC$ is spherical. Also, by the rules of graphical calculus, $d_{X} = d_{X^*}$.

\subsection{Multi-fusion Categories}
\label{subsec:multifusion}
Let $\CatC$ be a rigid monoidal category. $\CatC$ is called $\Complex$-linear if all $\Hom$ sets are finite dimensional vector spaces over $\Complex$, and the composition and tensor product of morphisms are $\Complex$-linear w.r.t each component. An object $X$ in a $\Complex$-linear category is called simple if $\End(X) = \Complex\, id_X$. An idempotent, i.e., a morphism $f \in \End(Y)$ such that $f^2 = f$, is called split if there are morphisms $g \in \Hom(Z, Y), h \in \Hom(Y, Z)$ for some $Z$, such that $hg = id_{Z}$ and $gh = f$. A $\Complex$-linear category is called semi-simple if it has direct sums, all idempotents split, and every object is isomorphic to a direct sum of simple objects \cite{muger2003subfactors}. There is a unique zero object denoted by $0$ in a semi-simple category.

\begin{remark}
\begin{itemize}
\item Another way to define a semi-simple category is to require a priori the category to be Abelian. This is equivalent to the current definition \cite{muger2003subfactors}. We avoid to use the terminology \lq Abelian category' since we will not deal with kernels and cokernels.
\item Given a $\Complex$-linear category $\CatC$, there is a canonical way to embed it as a full subcategory into a category $\DSC{\CatC}$ which has direct sums \cite{gabriel1997representations}. Roughly speaking, to define $\DSC{\CatC}$ one just formally introduces direct sums of objects of $\CatC$ as objects of $\DSC{\CatC}$ and defines the morphism spaces in the most natural way. There is also a standard way, called idempotent completion (or Karoubi envelope or Cauchy completion), to embed $\CatC$ into one $\IC{\CatC}$ which has all idempotents split. Moreover, $\IC{(\DSC{\CatC})} = \DSC{(\IC{\CatC})}$. We will discuss more about this in Section \ref{subsec:idempotent}.
\end{itemize}
\end{remark}

\begin{definition}
A multi-fusion category over $\Complex$ is a $\Complex$-linear rigid monoidal category which is semi-simple and has finitely many isomorphism classes of simple objects. A fusion category is a multi-fusion category in which the unit $\unit$ is simple.
\end{definition}

Let $\CatC$ be a multi-fusion category, and $\unit = \bigoplus\limits_{i \in I} \unit_i$ where the $\unit_i\,$'s are simple objects. One can show that $\unit_i \otimes \unit_j \simeq \delta_{i,j} \unit_i$. Moreover, for any simple object $X$, there is a unique $i \in I$, $j \in I$ such that $\unit_i \otimes X \simeq X \simeq X \otimes \unit_j$. Let $\CatC_{ij}$ be the full subcategory spanned by such simple objects. Then we have
\begin{align*}
\CatC = \bigoplus\limits_{i,j \in I} \CatC_{ij}.
\end{align*}
It follows that $\CatC_{ik} \otimes \CatC_{kj} \subset \CatC_{ij}$. Each $\CatC_{ii}$ is a fusion category with the unit $\unit_i$ and each $\CatC_{ij}$ is a $\CatC_{ii}$-$\CatC_{jj}$ bi-module category.
We call $\CatC$ an $|I| \times |I|$ multi-fusion category with index set $I$, $\CatC_{ij}$ the sector indexed by $(i,j)$, and call an object homogeneous if it belongs to some sector. If two homogeneous objects are from different sectors, then the only morphism between them is $0$. More generally, if $X = \bigoplus\limits_{i,j \in I}X_{ij}, \, Y = \bigoplus\limits_{i,j \in I}Y_{ij}, X_{ij}, Y_{ij} \in \CatC_{ij}$, then any morphism $f \in \Hom(X,Y)$ can be written as $f = (f_{ij})_{i,j \in I}$, where $f_{ij} \in \Hom(X_{ij}, Y_{ij})$.

Here are some examples of multi-fusion categories.
\begin{example}
\label{example}
\begin{enumerate}
\item \textbf{The \boldmath{$n \times n$}-matrix \boldmath{$\mathcal{M}_n$}}: the index set is $I = \{1,2,\cdots, n\}$. Each $(i,j)$-sector contains exactly one simple object $E_{ij}$. The tensor product obeys matrix multiplication rule: $E_{ik} \otimes E_{k'j} = \delta_{k,k'}E_{ij}$. Moreover, $\unit_i = E_{ii}, \, E_{ij}^* = E_{ji}$. All structural isomorphisms and $b_A\,'$s, $d_{A}\,'$s are the identity map.
\item \textbf{\boldmath{$G$}-graded fusion category}: For a finite group $G$, let $\CatC = \bigoplus\limits_{g} \CatC_g$ be a $G$-graded fusion category, that is, a fusion category such that $\CatC_{g} \otimes \CatC_{g'} \subset \CatC_{gg'}$. Define a multi-fusion category $\tilde{\CatC}$ whose index set is $G$ and whose $(g,g')$-sector is $\CatC_{g^{-1}g'}$. The tensor product and dual in $\tilde{\CatC}$ are the same as those in $\CatC$. Note that all the diagonals $\tilde{\CatC}_{gg}$ are copies of $\CatC_{e}$.
\end{enumerate}
\end{example}

Now let $\CatC$ be a multi-fusion category which is also pivotal. Note that $\unit$ is not a homogeneous object. If $A \in \CatC_{ij}$, then $A^* \in \CatC_{ji},\,  A \otimes A^* \in CatC_{ii}$, thus the birth map $b_A: \unit \longrightarrow A \otimes A^*$ can be equivalently viewed as a morphism in $\Hom(\unit_i, A \otimes A^*)$ since all other components of $b_A$ are zero. Similarly, we regard $b_{A}' \in \Hom(\unit_j, A^* \otimes A),\, d_A \in \Hom(A^* \otimes A, \unit_j),\, d_A' \in \Hom(A \otimes A^*, \unit_i)$.

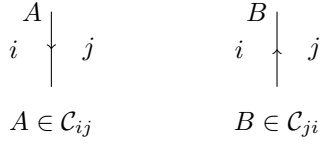
\begin{figure}
\centering
\begin{tikzpicture}[scale = 0.5]
 \begin{scope}[decoration={
    markings,
    mark=at position 0.5 with {\arrow{>}}}]
    \draw [postaction={decorate}] (0,2)node[left]{$A$} -- (0,0);
    \draw (-1,1)node{$i$};
    \draw (1,1)node{$j$};
    \draw (0,-1)node{$A \in \CatC_{ij}$};
 \end{scope}
 \begin{scope}[xshift = 6cm,decoration={
    markings,
    mark=at position 0.5 with {\arrow{>}}}]
    \draw [postaction={decorate}] (0,0) -- (0,2)node[left]{$B$};
    \draw (-1,1)node{$i$};
    \draw (1,1)node{$j$};
    \draw (0,-1)node{$B \in \CatC_{ji}$};
 \end{scope}
\end{tikzpicture}
\caption{Rules of labeling an arc in a graph diagram.}\label{fig:arc}
\end{figure}

Since $\CatC$ is pivotal, graphical calculus still makes sense in $\Real \times [0,1]$. However, since we will be mostly interested in homogeneous objects and also to avoid zero object, we refine the notion of graphical calculus reviewed in Section \ref{subsec:pivotal}. Let $\CatG$ be a graph diagram in $\Real \times [0,1]$ endowed with the standard counter clockwise orientation. The complement of $\CatG$ is divided into a disjoint union of connected regions, which are denoted by $R_1, R_2, \cdots$. Then a colored graph diagram is a graph diagram $\CatG$ satisfying the following condition:
\begin{itemize}
\item Each $R_i$ is labeled by an index $g_i \in I$, each arc is labeled by an object of $\CatC$, and each rectangle is labeled by a morphism obeying the following rules. For each arc $\beta$, let $R_i, R_j$ be the two regions bounded by $\beta$ such that the direction of $\beta$ together with the arrow pointing from $R_i$ to $R_j$ matches the orientation of $\Real \times [0,1]$. Then we require the object labeling $\beta$ to be from $\CatC_{g_i g_j}$. See Figure \ref{fig:arc} for an illustration. The requirement on the labeling of rectangles is the same as that mentioned in Section \ref{subsec:pivotal}.
\end{itemize}

The rules for interpreting colored graph diagrams as morphisms are the same as before. One can check the additional requirement on the labeling of arcs are consistent with these rules. For instance, if $A \in \CatC_{ij}$ labels an arc $\beta$, then reversing the direction of $\beta$ and changing $A$ to $A^* \in \CatC_{ji}$ still make it a well-defined coloring.

Let $X \in \CatC_{ij}, f \in \End(X)$, then it is direct to see that $\Tr^l(f) \in \End(\unit_j), \Tr^r(f) \in \End(\unit_i)$. See Figure \ref{fig:trace2}. Therefore, if $i \neq j$, then $\Tr^l(f)$ can never be equal to $\Tr^r(f)$. We conclude that a pivotal $n \times n$ multi-fusion category for $n > 1$ cannot be spherical according to the existing definition of sphericity. However, since the $\unit_i\,'$s are simple, we have $\Tr^l(f) = |\Tr^l(f)| id_{\unit_j}, \Tr^r(f) = |\Tr^r(f)| id_{\unit_i}$ for some complex numbers $|\Tr^l(f)|,\, |\Tr^r(f)|$. When interpreting a closed colored graph diagram as a morphism in some $\End(\unit_i)$ which is equal some complex number times $id_{\unit_i}$, what we are really interested in is the complex number but not the morphism itself. This motivates us to propose the following weakened definition:
\begin{definition}
\label{def:spherical_multi}
A spherical multi-fusion category $(\SMFC)$ is a pivotal multi-fusion category $\CatC$ such that $|\Tr^l(f)| = |\Tr^r(f)|$ for all $f \in \End(X), X \in \CatC_{ij}$. Define the trace of $f$ to be $\Tr(f):= |\Tr^l(f)|$, and the dimension of $X$ to be $d_{X}:= \Tr(id_{X})$.
\end{definition}
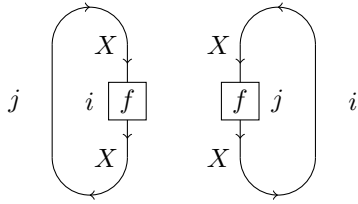
\begin{figure}
\centering
\begin{tikzpicture}[scale = 0.5]
 \begin{scope}[decoration={
    markings,
    mark=at position 0.5 with {\arrow{>}}}]
    \draw [postaction={decorate}] (2,1) arc(0:-180:1cm);
    \draw [postaction={decorate}] (0,4) arc(180:0:1cm);
    \begin{scope}[xshift = 2cm, yshift = 1cm]
      \draw [postaction={decorate}] (0,1) -- (0,0)node[left]{$X$};
      \draw (-0.5,1) rectangle (0.5, 2) node[pos = 0.5]{$f$};
      \draw [postaction={decorate}]  (0,3) node[left]{$X$} -- (0,2);
    \end{scope}
    \draw (0,1) -- (0,4);
    \draw (-1,2.5)node{$j$};
    \draw (1,2.5)node{$i$};
 \end{scope}

 \begin{scope}[xshift = 5cm, decoration={
    markings,
    mark=at position 0.5 with {\arrow{>}}}]
    \draw [postaction={decorate}] (0,1) arc(-180:0:1cm);
    \draw [postaction={decorate}] (2,4) arc(0:180:1cm);
    \begin{scope}[xshift = 0cm, yshift = 1cm]
      \draw [postaction={decorate}] (0,1) -- (0,0)node[left]{$X$};
      \draw (-0.5,1) rectangle (0.5, 2) node[pos = 0.5]{$f$};
      \draw [postaction={decorate}]  (0,3) node[left]{$X$} -- (0,2);
    \end{scope}
    \draw (2,1) -- (2,4);
    \draw (3,2.5)node{$i$};
    \draw (1,2.5)node{$j$};
 \end{scope}
 \end{tikzpicture}
 \caption{$\Tr^l(f)$ (Left) and $\Tr^r(f)$ (Right) for $X \in \CatC_{ij}$}\label{fig:trace2}
\end{figure}

For fusion categories, the above definition coincides with the existing definition of sphericity. Just as in the case of spherical categories, graphical calculus in a $\SMFC$ can also be generalized from the plane to $\Sphere^2$. One point to keep in mind is that when interpreting closed diagrams, it is the scalar  but not the morphism that remains invariant under isotopy.


At the end of this section, we introduce a special class of $\SMFC$s, which are the ingredients that will be used to construct invariants of $3$-manifolds in Section \ref{sec:main}. Let $\Label{\CatC}$ be a complete set of representatives, i.e., a set of simple objects that contains exactly one representative from each isomorphism class of simple objects, and let $\Label{\CatC}_{ij}$ be the subset of $\Label{\CatC}$ whose objects are from $\CatC_{ij}$, then $\Label{\CatC} = \bigsqcup\limits_{i,j \in I} \Label{\CatC}_{ij}$. Define the dimension of $\CatC_{ij}$ to be $K(\CatC_{ij}):= \sum\limits_{a \in \Label{\CatC}_{ij}} d_a^2$, the dimension of the $i$-th row to be $K(\CatC_i):= \sum\limits_{j \in I}K(\CatC_{ij})$, and the dimension of $\CatC$ to be $K(\CatC) = \sum\limits_{i,j \in I} K(\CatC_{ij})$.

\begin{definition}
A $\SMFC$ $\CatC$ is called $\Special$ if $K(\CatC_i)$ is the same for all $i \in I$.
\end{definition}


\section{Construction of Quantum Invariants}
\label{sec:main}
The Turaev-Viro-Barrett-Westbury $(\TVBW)$ invariant is a quantum invariant of $3$-manifolds constructed from a spherical fusion category. In this section, we generalize this construction to produce a quantum invariant of $3$-manifolds from a $\Special$ $\SMFC$ which is defined in Section \ref{subsec:multifusion}.
  
Let $\CatC = \bigoplus\limits_{i,j \in I} \CatC_{ij}$ be a $\Special$ $\SMFC$ with index set $I$. Recall that $\CatC$ is $\Special$ if $K(\CatC_i)$ is the same for all $i \in I$. In this case, denote $K(\CatC_i)$ by $K$. Let $\Label{\CatC}$ be a complete set of representatives, i.e., a set of simple objects that contains exactly one representative from each isomorphism class of simple objects, and let $\Label{\CatC}_{ij}$ be the subset of $\Label{\CatC}$ whose objects are from $\CatC_{ij}$, then $\Label{\CatC} = \bigsqcup\limits_{i,j \in I} \Label{\CatC}_{ij}$. By definition, $K(\CatC_i) = \sum\limits_{j \in I} \sum\limits_{a \in \Label{\CatC}_{ij}} d_a^2$. For any two homogeneous objects $X,Y \in \CatC_{ij}$, a pairing on $\Hom(X,Y) \times \Hom(Y,X)$ is defined as:
\begin{align}
\label{equ:pairingdef}
\langle\;,\;\rangle: \Hom(X,Y) \times \Hom(Y,X) &\longrightarrow   \Complex \nonumber \\
                       (\ \phi\quad\;  , \quad\psi\ ) \quad\qquad        &\mapsto           \Tr(\phi\psi)
\end{align}
Recall from Section \ref{subsec:multifusion} that $\Tr(\phi\psi) = |\Tr^l(\phi\psi)| = |\Tr^r(\phi\psi)|$. The pairing is non-degenerate. Thus, there are natural isomorphisms $\Hom(X,Y)\simeq \Hom(Y,X)^*$, $\Hom(Y,X) \simeq \Hom(X,Y)^*$.

The $\TVBW$ invariant from spherical fusion categories is defined on a triangulation of $3$-manifolds \cite{barrett1996invariants}, or more generally on a polytope decomposition \cite{kirillov2010turaev}. Here the invariant to be introduced below can also be defined both on triangulations and polytope decompositions. For simplicity, we restrict ourselves on triangulations.

Let $M$ be a closed oriented $3$-manifolds. By a triangulation of $M$ is meant a $\Delta$-complex whose underlying space is homeomorphic to $M$. An ordered triangulation is one whose vertices are ordered by $0,1, \cdots$. Let $\CatT$ be an ordered triangulation of $M$. Denote by $\CatT^i$ the set of $i$-simplices of $\CatT$. For each $i$-simplex $\sigma$ of $\CatT$, the ordering on $\CatT^0$ induces a relative ordering on the vertices of $\sigma$ with which we can identify $\sigma$ with the standard $i$-simplex $(0,1,\cdots, i)$. The invariant of $3$-manifolds to be defined will only depend on this relative ordering for each simplex.

\begin{definition}
Let $M, \CatT, \CatC$ be as above and $\Label{\CatC}$ be an arbitrary complete set of representatives. A $\CatC$-coloring of the pair $(M, \CatT)$ is a pair of functions $F = (f^0, f^1)$, $f^0: \CatT^0 \longrightarrow I$, $f^1: \CatT^1 \longrightarrow \Label{\CatC}$, such that for each $1$-simplex $\alpha = (01)$ with the induced ordering on its vertices,
\begin{align*}
f^1_{01} \in \Label{\CatC}_{f^0_0,f^0_1}.
\end{align*}
\end{definition}
In the above definition, we have identified a $1$-simplex with the standard one $(01)$. Under the absolute ordering, a $1$-simplex whose vertices are ordered by $(i,j), \ i < j$ is subject to the condition that $f^1_{ij} \in \Label{\CatC}_{f^0_i,f^0_j}$. In the following we will use this identification for other simplices as well.

Assume a coloring $F$ has been given. For each $2$-simplex $\beta = (012)$, we have $f^1_{ij} \in \Label{\CatC}_{f^0_i,f^0_j}$, $0 \leq i < j \leq 2$, then $f^1_{01} \otimes f^1_{12} \in \Label{\CatC}_{f^0_0,f^0_2}$ is in the same sector as $f^1_{02}$. Define,
\begin{align*}
V^{+}_{F}(\beta) = \Hom(f^1_{02}, f^1_{01} \otimes f^1_{12}), \qquad, V^{-}_{F}(\beta) = \Hom(f^1_{01} \otimes f^1_{12},f^1_{02}).
\end{align*}
Then by the non-degenerate pairing in Equation \ref{equ:pairingdef}, $V^{+}_{F}(\beta) \simeq V^{-}_{F}(\beta)^*$, $V^{-}_{F}(\beta) \simeq V^{+}_{F}(\beta)^*$.

\begin{figure}
\centering
\includegraphics[scale=1]{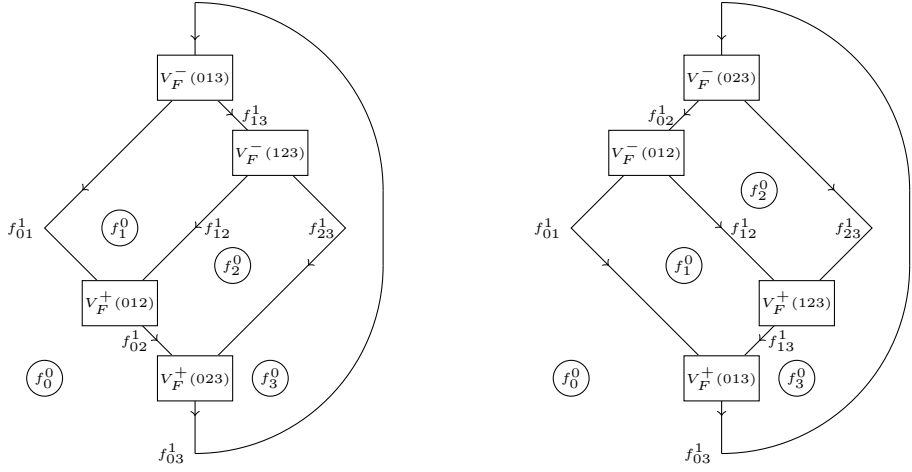}
\caption{Graphical definition of $\tilde{Z}^{+}_{F}(0123)$ (Left) and $\tilde{Z}^{-}_{F}(0123)$ (Right). Here we put the vertex colors in circles to distinguish with edge colors.}\label{fig:6jbox}
\end{figure}

For each $3$-simplex $\gamma = (0123)$, define a linear functional $\tilde{Z}^{+}_{F}(\gamma)$,
\begin{align*}
\tilde{Z}^{+}_{F}(\gamma): V^{-}_F(123) \otimes V^{+}_F(023) \otimes V^{-}_F(013) \otimes V^{+}_F(012) \longrightarrow \Complex
\end{align*}
as follows. For $\phi_{123} \otimes \phi_{023}\otimes \phi_{013}\otimes \phi_{012}$ in the domain, define,
\begin{align}
\label{equ:Zptilde}
\tilde{Z}^{+}_{F}(\gamma) (\phi_{123} \otimes \phi_{023}\otimes \phi_{013}\otimes \phi_{012}) = \Tr(\phi_{012}(id \otimes \phi_{123})(\phi_{012} \otimes id)\phi_{023}),
\end{align}
or graphically as the diagram given in Figure \ref{fig:6jbox} (Left), where the value of $\tilde{Z}^{+}_{F}(\gamma)$ is the evaluation of the diagram with each box colored by the corresponding $\phi_{ijk}$. One can check the requirement on the coloring makes the composition of the $\phi_{ijk}\,$'s in Equation \ref{equ:Zptilde} well defined. By the non-degenerate pairing, $\tilde{Z}^{+}_{F}(\gamma)$ induces a linear map,
\begin{align*}
Z^{+}_{F}(\gamma): V^{+}_F(023) \otimes V^{+}_F(012) \longrightarrow V^{+}_F(123) \otimes V^{+}_F(013),
\end{align*}
such that
\begin{align*}
\langle Z^{+}_{F}(\gamma)(\phi_{023} \otimes \phi_{012}), \phi_{123} \otimes \phi_{013} \rangle = \tilde{Z}^{+}_{F}(\gamma) (\phi_{123} \otimes \phi_{023}\otimes \phi_{013}\otimes \phi_{012}).
\end{align*}
Similarly, we define
\begin{align*}
\tilde{Z}^{-}_{F}(\gamma): V^{+}_F(123) \otimes V^{-}_F(023) \otimes V^{+}_F(013) \otimes V^{-}_F(012) \longrightarrow \Complex
\end{align*}
by
\begin{align}
\label{equ:Zmtilde}
\tilde{Z}^{-}_{F}(\gamma) (\phi_{123}' \otimes \phi_{023}'\otimes \phi_{013}'\otimes \phi_{012}') = \Tr(\phi_{023}'(\phi_{012}' \otimes id)(id \otimes \phi_{123}')\phi_{013}'),
\end{align}
or graphically as the diagram shown in Figure \ref{fig:6jbox} (Right). And in the same way, this induces a linear map,
\begin{align*}
Z^{-}_{F}(\gamma): V^{+}_F(123) \otimes V^{+}_F(013) \longrightarrow V^{+}_F(023) \otimes V^{+}_F(012).
\end{align*}

Let $V^{-}_F(\gamma) = V^{+}_F(023) \otimes V^{+}_F(012)$ and $V^{+}_F(\gamma) = V^{+}_F(123) \otimes V^{+}_F(013)$. Define $\epsilon(\gamma) = +$ if the orientation on $\gamma$ induced from that of $M$ coincides with the one determined by the ordering of its vertices, and define $\epsilon(\gamma) = -$ otherwise. Then we have
\begin{align*}
Z^{\epsilon(\gamma)}_{F}(\gamma): V^{-\epsilon(\gamma)}_F(\gamma) \longrightarrow V^{+\epsilon(\gamma)}_F(\gamma),
\end{align*}
where we adopt the convention $++ = -- = +,\, +- = -+ = -$.

One observation on the definition of $\epsilon(\gamma)$ is as follows. If $\gamma = (0123)$ is a $3$-simplex, then $\epsilon(\gamma)(0123)$ matches the orientation of $M$. A boundary face $\beta = (ijk)$ of $\gamma$ is called positive if its orientation induced by the ordering of its vertices matches $\partial (\epsilon(\gamma)(0123))$, and is called negative otherwise. Then $V_F^{+}(\beta)$ appears as a component in the domain of $Z_F^{\epsilon(\gamma)}(\gamma)$ if it is negative and as a component in the codomain otherwise.

Let $V_F = \bigotimes\limits_{\beta \in \CatT^2} V^{+}_F(\beta)$. Since $M$ is closed, each $2$-simplex $\beta$ is the common face of exactly two $3$-simplices $\gamma_1, \gamma_2$ ($\gamma_1$ could be the same as $\gamma_2$), and moreover, the sign of $\beta$ in $\gamma_1$ and $\gamma_2$ are opposite. Thus, if $V^{+}_F(\beta)$ appears as a component in the domain of $Z^{\epsilon(\gamma_1)}_{F}(\gamma_1)$, then it must appear as a component in the codomain of $Z^{\epsilon(\gamma_2)}_{F}(\gamma_2)$. From this observation, we have $V_F = \bigotimes\limits_{\gamma \in \CatT^3} V^{-\epsilon(\gamma)}_F(\gamma) =\bigotimes\limits_{\gamma \in \CatT^3} V^{+\epsilon(\gamma)}_F(\gamma)$ (up to permutation of tensor components or viewed as an unordered tensor product). This implies $\bigotimes\limits_{\gamma \in \CatT^3} Z^{\epsilon(\gamma)}_{F}(\gamma)$ is an endmorphism on $V_F$.

\begin{definition}
Let $M, \CatT, \CatC$ be as above. The partition function of the pair $(M, \CatT)$ is defined to be,
\begin{align}
\label{equ:partition}
Z_{\CatC}(M,\CatT) = \sum\limits_{F = (f^0,f^1)} K^{-|\CatT^0|} \prod\limits_{\alpha \in \CatT^1} d_{f^1_{\alpha}} \,  \Tr(\bigotimes\limits_{\gamma \in \CatT^3} Z^{\epsilon(\gamma)}_{F}(\gamma)),
\end{align}
where the summation is over all colorings.
\end{definition}

\begin{remark}
In the $\TVBW$ construction, the relevant factor involving $|\CatT^0|$ is $K(\CatC')^{-|\CatT^0|}$, where $K(\CatC')$ is the dimension of a spherical fusion category $\CatC'$. However, in the definition of the current invariant, the corresponding factor is $K^{-|\CatT^0|}$, where $K = K(\CatC_i)$ is the dimension of the $i$-th row, i.e., the direct sum $\bigoplus\limits_{j\in I}\CatC_{ij}$, and we require $K$ to be independent of $i \in I$. This requirement is necessary when in the proof of invariance of the partition function under the Pachner $1$-$4$ move.
\end{remark}

The main result is as follows.
\begin{theorem}
\label{thm:main}
The partition function $Z_{\CatC}(M,\CatT)$ is independent of the choice of $a)$ the complete set of representatives $\Label{\CatC}$, $b)$ the ordering on the vertices of $\CatT$, and $c)$ the triangulation $\CatT$. Therefore, $Z_{\CatC}(M):= Z_{\CatC}(M,\CatT)$ is an invariant of closed oriented $3$-manifolds.
\begin{proof}
The proof is mostly parallel to that in the case of spherical fusion categories in \cite{barrett1996invariants}. To avoid repetition, we only illustrate why $\CatC$ is required to be $\Special$. Let $(01234)$ be the standard $4$-simplex whose boundary is partitioned into $\CatT \sqcup \CatT'$ with $\CatT = (0124) \sqcup (0234), \CatT' = (1234) \cup (0134) \cup (0123) $. Then $\CatT$ and $\CatT'$ share all vertices and edges except that $\CatT'$ has one edge $(13)$ of its own. Let $F$ be a coloring on $\CatT$, and $F'$ be an extension from $F$ to a coloring on $\CatT'$. To emphasize the change of coloring on $(13)$, in the following we will write $\sum_{F'}$ as $\sum_{13}$, with the understanding that the colors on all other simplices are fixed. We also drop the subscript $F$ in $Z^{\pm}_{F}(\cdot)$. As in \cite{barrett1996invariants}, the following two facts hold.
\begin{align}
\label{equ:Z2-3}
Z^+(0124)_{2,3}Z^+(0234)_{1,2} = \sum\limits_{13}d_{13} Z^+(1234)_{1,2}Z^+(0134)_{1,3}Z^+(0123)_{2,3},
\end{align}
\begin{align}
\label{equ:Zorthogonal}
d_{02}\sum\limits_{13} d_{13} Z^{-}(0123)Z^{+}(0123) = id,
\end{align}
where in the first equation the map on either side is from $V^{+}(034) \otimes V^{+}(023)\otimes V^{+}(012)$ to $V^{+}(234) \otimes V^{+}(124)\otimes V^{+}(014)$, and $Z^+(\cdot)_{i,j}$ means $Z^+(\cdot)$ acts on the $i$-th and $j$-th components.

The invariance of $Z_{\CatC}(M)$ under pachner move $2$-$3$ is proved with Equation \ref{equ:Z2-3}. To prove the invariance under pachner move $1$-$4$, it suffices to show,
\begin{align}
\label{equ:Z1-4}
& Z^+(0234) \\
 =& \sum\limits_{1,01,12,13,14}K^{-1} d_{01}d_{12}d_{13}d_{14}  \Tr_3(Z^-(0124)_{2,3}Z^+(1234)_{1,2}Z^+(0134)_{1,3}Z^+(0123)_{2,3}) \nonumber,
\end{align}
where $\Tr_3(\cdot)$ means taking partial trace with respect to the $3$rd component, and the summation on the right side is over all colorings which change colors only on the vertex $1$ and on the edges $(01),(12),(13),(14)$. We prove Equation \ref{equ:Z1-4}.
\begin{align*}
\text{RHS} &\overset{Equ. \ref{equ:Z2-3}}{=} \sum\limits_{1,01,12,14}K^{-1} d_{01}d_{12}d_{14}  \Tr_3(Z^-(0124)_{2,3}Z^+(0124)_{2,3}Z^+(0234)_{1,2}) \\
           &\overset{Equ. \ref{equ:Zorthogonal}}{=} \sum\limits_{1,01,12}K^{-1} d_{01}d_{12}d_{02}^{-1}  \Tr_3(Z^+(0234)_{1,2}) \\
           &=\sum\limits_{1,01,12}K^{-1} d_{01}d_{12}d_{02}^{-1}\dim {V^+(012)} \quad Z^+(0234)\\
           &= Z^+(0234) \\
\end{align*}
The last equality above is due to the following property. For fixed $c \in \Label{\CatC}_{ij}$,
\begin{align}
\label{equ:fixedc}
\sum\limits_{k \in I}\sum\limits_{\substack{a \in \Label{\CatC}_{ik}, \\ b \in \Label{\CatC}_{kj}}} d_ad_bd_c^{-1} N_{ab}^c =& \sum\limits_{k \in I}\sum\limits_{\substack{a \in \Label{\CatC}_{ik}, \\ b \in \Label{\CatC}_{kj}}} d_ad_bd_{\bar{c}}^{-1} N_{\bar{c}a}^{\bar{b}} =& \sum\limits_{k \in I}\sum\limits_{a \in \Label{\CatC}_{ik}} d_a^2 &= K(\CatC_i),
\end{align}
where $N_{ab}^c:= \dim \Hom(a\otimes b,c),\, \bar{c}:= c^{*}$.

\end{proof}
\end{theorem}

\begin{example}
As an example, we compute the invariant of $\Sphere^3$ with the standard orientation. We use the notations in the proof of Theorem \ref{thm:main}. Let $\gamma_1 = (0123), \gamma_2 = (0123)$ be the two standard $3$-simplices glued together along their corresponding faces, one with positive orientation and the other with negative orientation. Thus their union $\gamma_1 \cup \gamma_2$ is a triangulation of $\Sphere^3$. With this triangulation we have, 
\begin{align*}
Z_{\CatC}(\Sphere^3) &= \sum\limits_{\substack{0,1,2,3 \\ 01,02,03,12,13,23}} K^{-4} \prod\limits_{(i,j):0 \leq i<j\leq 3} d_{ij} \quad \Tr(Z^{-}(0123)Z^{+}(0123)) \\
                     &\overset{Equ. \ref{equ:Zorthogonal}}{=} \sum\limits_{\substack{0,1,2,3 \\ 01,02,03,12,23}} K^{-4} d_{01}d_{03}d_{12}d_{23} \dim V^+(023) \dim V^+(012) \\
                     &\overset{Equ. \ref{equ:fixedc}}{=} \sum\limits_{\substack{0,2,3 \\02,03,23}} K^{-3} d_{03}d_{23}d_{02} \dim V^+(023) \\
                     &\overset{Equ. \ref{equ:fixedc}}{=} \sum\limits_{\substack{0,3 \\03}} K^{-2} d_{03}^2 \\
                     &= \frac{|I|}{K}.
\end{align*}
\end{example}

\begin{figure}
\centering
\includegraphics[scale=1]{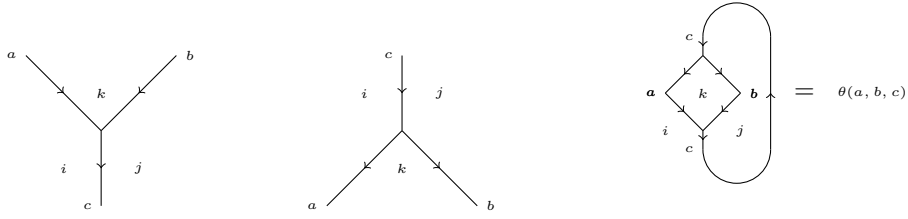}
\caption{Graphical representation of $B_{c}^{ab}$ (Left), $B_{ab}^c$ (Middle), and their relation (Right). Here $i,j,k$ are indices coloring regions}\label{fig:basis}
\end{figure}

In the following, we give a formula for $Z_F = \Tr(\bigotimes\limits_{\gamma \in \CatT^3} Z^{\epsilon(\gamma)}_{F}(\gamma))$ under certain basis. For simplicity, we assume the category $\CatC$ is multiplicity free, that is, for any three simple objects $a,b,c$, $\Hom(c, a \otimes b)$ has dimension either $0$ or $1$. If it is the latter case, we call $(a,b,c)$ {\it admissible}. Now, for any admissible $(a,b,c)$, we choose a basis element $B_{c}^{ab} \in  \Hom(c, a\otimes b)$ and $B_{ab}^c \in \Hom(a \otimes b, c)$ such that,
\begin{align*}
\langle B_{c}^{ab} , B_{ab}^c \rangle = \Tr(B_{ab}^cB_{c}^{ab}) = \theta(a,b,c),
\end{align*}
where $\theta(a,b,c)$ is certain constant to be specified later. Graphically, $B_{c}^{ab}, B_{ab}^c$ and their relations are represented as in Figure \ref{fig:basis}.

\begin{figure}
\centering
\includegraphics[scale=1]{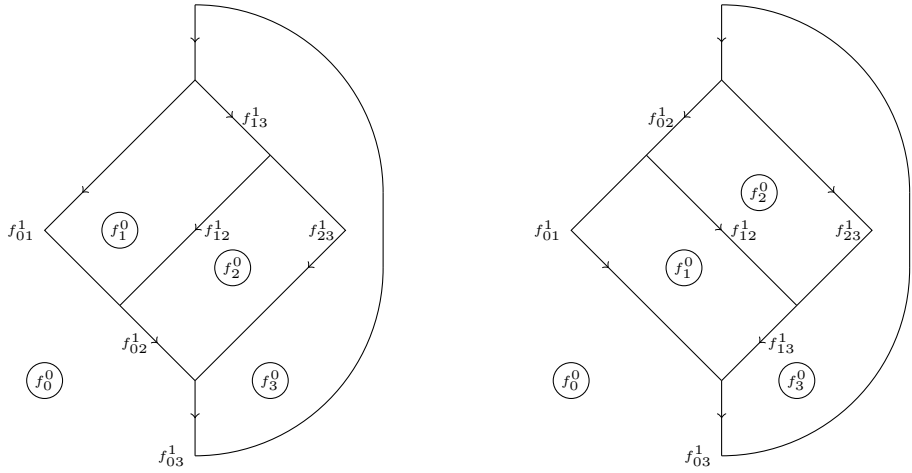}
\caption{Graphical representations of $\tilde{Z}^{+}_F(0123;B)$ (Left) and $\tilde{Z}^{-}_F(0123;B)$ (Right). Here we put the vertex colors in circles to distinguish with edge colors.}\label{fig:6j}
\end{figure}

\begin{proposition}
With the notations as above, the invariant has the \lq state-sum' formula:
\begin{align}
\label{equ:partition2}
Z_{\CatC}(M,\CatT) = \sum\limits_{F = (f^0,f^1)} \prod\limits_{\tau \in \CatT^0} K^{-1} \prod\limits_{\alpha \in \CatT^1} d_{f^1_{\alpha}} \,  \prod\limits_{\beta \in \CatT^2} \theta_F(\beta)^{-1} \, \prod\limits_{\gamma \in \CatT^3} \tilde{Z}^{\epsilon(\gamma)}_F(\gamma;B),
\end{align}
where $\tilde{Z}^{\epsilon(\gamma)}_F(\gamma;B)$ is defined as the evaluation of the diagrams in Figure \ref{fig:6j}.
\begin{proof}
Fix a coloring $F$ on the triangulation $\CatT$. For each $2$-simplex $\beta = (012)$, denote by,
\begin{align*}
B^{+}(\beta) = B_{f^1_{02}}^{f^1_{01}f^1_{12}} \in V^{+}_{F}(\beta), \quad B^{-}(\beta) = B^{f^1_{02}}_{f^1_{01}f^1_{12}} \in V^{-}_{F}(\beta), \quad \theta_F(\beta) = \langle B^+(\beta),B^-(\beta)\rangle.
\end{align*}
Then it follows that, for any $3$-simplex $\gamma = (0123)$,
\begin{align*}
Z_F^{+}(0123)(B^+(023)\otimes B^+(012)) = \frac{\tilde{Z}^+_F(0123;B)}{\theta_F(123)\theta_F(013)} B^+(123)\otimes B^+(013),
\end{align*}
where $\tilde{Z}^+_F(0123;B) = \tilde{Z}_F^{+}(0123)(B^-(123) \otimes B^+(023)\otimes B^-(013) \times B^+(012))$ is the evaluation of the colored graph in Figure \ref{fig:6j} (Left). Similarly, we have
\begin{align*}
Z_F^{-}(0123)(B^+(123)\otimes B^+(013)) = \frac{\tilde{Z}^-_F(0123;B)}{\theta_F(023)\theta_F(012)} B^+(023)\otimes B^+(012),
\end{align*}
where $\tilde{Z}^-_F(0123;B) = \tilde{Z}_F^{-}(0123)(B^+(123) \otimes B^-(023)\otimes B^+(013) \times B^-(012))$ is the evaluation of the colored graph in Figure \ref{fig:6j} (Right).
Then we have
\begin{align*}
Z_F = \Tr(\bigotimes\limits_{\gamma \in \CatT^3} Z^{\epsilon(\gamma)}_{F}(\gamma)) =  \prod\limits_{\beta \in \CatT^2} \theta_F(\beta)^{-1} \, \prod\limits_{\gamma \in \CatT^3} \tilde{Z}^{\epsilon(\gamma)}_F(\gamma;B).
\end{align*}
\end{proof}
\end{proposition}

A common choice of $\theta(a,b,c)$ is to have $\theta(a,b,c) = 1$ for all admissible $(a,b,c)$, in which case the formula in Equation \ref{equ:partition2} does not involve contributions from $2$-simplices. Another common choice in physics literature, assuming $\CatC$ is unitary which implies the quantum dimension of any non-zero object is positive, is to have $\theta(a,b,c) = \sqrt{d_ad_bd_c}$.

\section{Invariants from Generalized Categorical Groups}
\label{sec:generalized}
In this section, we study a class of $\Special$ $\SMFC$s obtained from what we call generalized categorical groups. We first have a review of $2$-groups and categorical groups, and then introduce the notion of generalized categorical groups which are generalizations of categorical groups. By a process called idempotent completion, we turn a generalized categorical group into a $\Special$ $\SMFC$. The partition functions (and $(2+1)$-$\TQFT$s) from such $\SMFC$s are shown to contain the ones in \cite{kapustin2013higher}, while the latter $\TQFT$ was not known to have a categorical construction before.  

\subsection{2-groups}
\label{subsec:2group}
A $2$-group is a triple $\CatG = (G, A, \beta)$, where $G$ is a finite group, $A$ is a finite Abelian group endowed with a $G$-action, and $\beta \in H^3(G, A)$ is a $3$rd cohomology class, where the cohomology group $H^3(G, A)$ is defined with respect to the $G$-action on $A$. By abusing of languages, we also assume $\beta$ is a co-cycle representing the class $\beta$. Different choices of representative co-cycles correspond to equivalent $2$-groups. We write the product in $A$ multiplicatively instead of additively. Denote the unit element in a group by $e$ and the inverse of an element $g$ by $g^{-1}$ or $\overline{g}$. Then $\beta$ being a co-cycle means that for any $g_0, g_1, g_2, g_3 \in G$,
\begin{align}
\label{equ:3rd cocycle}
e &= \delta \beta (g_1,g_2,g_3,g_4) \nonumber \\
&= \lsup{g_0}{\beta(g_1, g_2, g_3)} \overline{\beta(g_0g_1, g_2, g_3)} \beta(g_0, g_1g_2, g_3) \overline{\beta(g_0, g_1, g_2g_3)} \beta(g_0, g_1, g_2) .
\end{align}

A typical example of a $2$-group comes from the homotopy $2$-type $(\pi_1(X), \pi_2(X), \beta)$ of a complex $X$, where $\pi_2(X)$ is endowed with the monodromy action of $\pi_1(X)$ and $\beta$ is the Postnikov invariant of $X$ \cite{maclane1949cohomology}. Actually, $2$-groups classify homotopy $2$-types \cite{maclane19503}.

A categorical group is a rigid monoidal category in which all objects and all morphisms are invertible. From a $2$-group $\CatG = (G, A, \beta)$, a categorical group $\CatC(\CatG)$ can be constructed as follows.
\begin{enumerate}
\item $\CatC(\CatG)^0 = G$; $\Hom_{\CatC(\CatG)}(g_1, g_2)$ is $A$ if $g_1 = g_2$ and the empty set otherwise; composition of morhisms is multiplication in $A$.
\item For $g_1, g_2 \in G, h_1 \in \Hom(g_1, g_1), h_2 \in \Hom(g_2, g_2)$, $g_1 \otimes g_2:= g_1g_2$, $h_1 \otimes h_2:= h_1 \lsup{g_1}{(h_2)}$, $\unit:= e$.
\item For $g_1, g_2, g_3 \in G,$ the association isomorphism is defined to be
\begin{equation*}
\beta(g_1,g_2,g_3): (g_1 \otimes g_2) \otimes g_3 \longrightarrow g_1 \otimes (g_2 \otimes g_3).
\end{equation*}
\item The dual of an object $g$ is $g^* := \overline{g}$.
\end{enumerate}

It is straight forward to check the above defines a categorical group. For instance, the Pentagon Equation that the association isomorphism needs to satisfy translates exactly to the co-cycle condition in Equation \ref{equ:3rd cocycle}. Actually it is also true that up to some appropriately defined equivalence, $2$-groups are in one-to-one correspondence with categorical groups \cite{Etingof_finitetensor}.

One can also \lq linearize' the categorical group $\CatC(\CatG)$ by redefining the objects as formal direct sums of elements of $G$ and the morphisms from $g \in G$ to itself as linear spans of elements of $A$ with coefficients in $\Complex$, namely, $\Hom(g,g) = \Complex[A]$. (If $g_1 \neq g_2$, then $\Hom(g_1,g_2)$ is redefined to be the zero vector space.) The composition and tensor product are extended linearly. We still denote the \lq linearized' category as $\CatC(\CatG)$. Note that if $A$ is not the trivial group, then $\CatC(\CatG)$ is not a semisimple category since $\Hom(g,g)$ is not isomorphic to $\Complex$. A more fundamental reason is that there are idempotents in $\Hom(g,g)$ which do not split. We show in Section \ref{subsec:idempotent} that by a process called idempotent completion, $\CatC(\CatG)$ can be turned into a semisimple category, or more specifically a $\SMFC$. Before doing that, we first show in Section \ref{subsec:generalized} that the notion of categorical groups can be generalized so that the tensor product and Pentagon solution encode more data that a three co-cycle.

\subsection{Generalized Categorical Groups}
\label{subsec:generalized}
Here we explore more general structures in a rigid monoidal category whose underlying objects and morphisms are the same as a categorical group. To be more precise, let $\CatC(G, A)$ be a rigid monoidal category such that the objects form a finite group $G$ by tensor product, $\Hom(g_1, g_2) = \delta_{g_1,g_2}\Complex[A]$ for $g_1, g_2 \in G$, and the composition of morphisms is multiplication in $\Complex[A]$, where $A$ is a finite Abelian group on which $G$ acts. Of-course, a categorical group arising from a $2$-group is such a category. We show below that a more general form of tensor product of morphisms and the solution to Pentagon Equation can be defined in $\CatC(G, A)$ beyond those in a categorical group.

Let $\hat{A}$ be the group complex characters on $A$. The action of $G$ on $A$ induces an action on $\hat{A}$ and extends by linearity to an action on $\Complex[A]$. Specifically, for $g \in G, \chi \in \hat{A}$, $\lsup{g}{\chi}: = \chi(\lsup{\overline{g}}{(\cdot)})$.

Let $h_1, h_2 \in A$ but consider $h_1 \in \Hom(g_1, g_1), h_2 \in \Hom(g_2,g_2)$, and define $h_1 \otimes h_2 \in \Hom(g_1g_2, g_1g_2)$ by
\begin{align*}
h_1 \otimes h_2:= \lambda(\overline{g_1})(h_2) h_1 \lsup{g_1}{h_2},
\end{align*}
where $\lambda: G \longrightarrow \hat{A}$ is some map. Extend the definition linearly to define tensor product of morphisms which are linear combinations of group elements.  Thus, compared to the tensor product in $\CatC(\CatG)$, an extra coefficient $\lambda(\overline{g_1})(h_2)$ is introduced in the current setting.

It is direct to check that the associativity $(h_1 \otimes h_2) \otimes h_3 = h_1 \otimes (h_2 \otimes h_3)$ is equivalent to the condition
\begin{align*}
\lambda(\overline{g_2} \, \overline{g_1}) = \lambda(\overline{g_2}) \lsup{\overline{g_2}}{\lambda(\overline{g_1})},
\end{align*}
 which means $\lambda$ is a co-cycle in $Z^1(G, \hat{A})$. Again, equivalent choices of $\lambda$ within the same cohomology would correspond to equivalent structures on the category, thus we can view $\lambda \in H^1(G, \hat{A})$. The case of $\CatC(\CatG)$ corresponds to the trivial cohomology class.

Before looking at the association isomorphisms, we first recall some properties of characters. For each $\chi \in \tilde{A}$, define $P_{\chi} \in \Complex[A]$ by
\begin{equation}
P_{\chi}:= \frac{1}{|A|} \sum\limits_{h \in A} \overline{\chi(h)}\, h.
\end{equation}
By standard character theories, the following properties hold.
\begin{itemize}
 \item $P_{\chi_1}P_{\chi_2} = \delta_{\chi_1, \chi_2} P_{\chi_1}$; $h = \sum\limits_{\chi \in \tilde{A}} \chi(h)P_{\chi}$; in particular, $e = \sum\limits_{\chi \in \tilde{A}} P_{\chi}$.
 \item $\lsup{g}{P_{\chi}} = P_{\lsup{g}{\chi}}$; $id_g \otimes P_{\chi} = P_{\lambda(g)\lsup{g}{\chi}}$.
\end{itemize}
The first property means $\{P_{\chi}: \chi \in \tilde{A}\}$ forms a set of complete orthogonal idempotents, and a basis of $\Complex[A]$ in particular.

A general element in $\Complex[A]$ is of the form $\sum\limits_{\chi}c_{\chi} P_{\chi}$, $c_{\chi} \in \Complex$, and it is invertible if and only if $c_{\chi} \neq 0$ for all $\chi$. Thus, for $g_1, g_2, g_3 \in G$, the association isomorphism $a(g_1,g_2,g_3)$ takes the form
\begin{align*}
a(g_1,g_2,g_3) = \sum\limits_{\chi} a(g_1,g_2,g_3)_{\chi} P_{\chi},\quad a(g_1,g_2,g_3)_{\chi} \neq 0.
\end{align*}
The Pentagon Equation is then equivalent to the condition,
\begin{align}
\label{equ:general_pentagon}
a(g_1g_2, g_3,g_4)_{\chi} a(g_1,g_2, g_3g_4)_{\chi} = a(g_1,g_2, g_3)_{\chi}a(g_1,g_2g_3,g_4)_{\chi}a(g_2, g_3,g_4)_{\lambda(\overline{g_1})\lsup{\overline{g_1}}{\chi}}.
\end{align}

Let $C^1(\hat{A})$ be the group of all maps from $\hat{A}$ to non-zero complex numbers $\Complex^{\times}$. Given $\lambda \in Z^1(G, \hat{A})$, we define a new action of $G$ on $\hat{A}$ by $\phi(g, \chi) := \lambda(g)\lsup{g}{\chi}$. Note that in this action, $\hat{A}$ is viewed as a set but not a group, and $\phi(g, \cdot)$ is not a group automorphism unless $\lambda$ is the trivial co-cycle. However, the induced action on $C^1(\hat{A})$ defined by $\phi(g, \psi) = \psi(\overline{g}, \cdot)$ is an action by automorphism whether or not $\lambda$ is trivial. We denote by $C^1(\hat{A})_{\phi}$ the group $C^1(\hat{A})$ with the induced action defined above. Then Equation \ref{equ:general_pentagon} is equivalent to
\begin{align}
a(g_1g_2, g_3,g_4) a(g_1,g_2, g_3g_4) = a(g_1,g_2, g_3)a(g_1,g_2g_3,g_4)\phi(g_1,a(g_2, g_3,g_4)),
\end{align}
where $a(g_1, g_2, g_3)$ is viewed as an element in $C^1(\hat{A})_{\phi}$. Hence, we have $a \in Z^3(G, C^1(\hat{A})_{\phi})$.

The above defined category depends on the data $(G, A, \lambda, a)$, where $\lambda \in H^1(G, \hat{A}), a \in H^3(G, C^1(\hat{A})_{\phi})$. We call such a category a {\it generalized categorical group} and denote it by $\CatC(G, A, \lambda, a)$.

If $\lambda$ is the trivial co-cycle, then the action $\phi: G \times \hat{A} \longrightarrow \hat{A}$ coincides with action of $G$ on $\hat{A}$ induced by the given action of $G$ on $A$, namely, $\phi(g, \chi) = \lsup{g}{\chi} = \chi(\lsup{\overline{g}}(\cdot))$. Note that $A \simeq \hat{\hat{A}}$ is a subgroup of $C^1(\hat{A})_{\phi}$. That is, given $h \in A, \chi \in \hat{A}$, $h$ is viewed as an element of $C^1(\hat{A})_{\phi}$ by $h(\chi) := \chi(h)$. If $\lambda$ is trivial, it can be shown that the embedding $\iota: A \hookrightarrow C^1(\hat{A})_{\phi}$ is $G$-equivariant. Then we have the induced map $\iota_*: H^3(G, A) \longrightarrow H^3(G, C^1(\hat{A})_{\phi})$. Let $a = \iota_*(\beta)$ for some $\beta \in H^3(G, A)$, then $a(g_1, g_2, g_3)_{\chi} = \chi(\beta(g_1,g_2,g_3))$, and hence
\begin{align*}
a(g_1, g_2, g_3) &= \sum\limits_{\chi} a(g_1, g_2, g_3)_{\chi} P_{\chi} &= \sum\limits_{\chi} \chi(\beta(g_1,g_2,g_3)) P_{\chi} &= \beta(g_1, g_2,g_3)
\end{align*}
Therefore, we recovered the categorical group constructed from the $2$-group $(G, A, \beta)$ when $\lambda$ is trivial and $a = \iota_*(\beta)$.

More generally if $\lambda$ is not necessarily trivial, let $\omega \in H^3(G,\Complex^{\times} ) \simeq H^3(G, U(1))$, $\beta \in H^3(G, A)$, and let
\begin{align}
\label{equ:a_omega_beta}
a(g_1,g_2,g_3)_{\chi}:= \omega(g_1,g_2,g_3) \chi(\beta(g_1,g_2,g_3)).
\end{align}
Then Equation \ref{equ:general_pentagon} can be rewritten as,
\begin{align}
\delta \omega(g_1,g_2,g_3,g_4) \cdot \chi (\delta\beta (g_1,g_2,g_3,g_4)) \cdot \lambda(\overline{g_1})(\beta(g_2,g_3,g_4)) = 1,
\end{align}
which is equivalent to
\begin{align}
\label{equ:cup}
\lambda(\overline{g_1})(\beta(g_2,g_3,g_4)) =\overline{\lambda(g_1)}(\lsup{g_1}{\beta(g_2,g_3,g_4)})  = 1,
\end{align}
where the first equality above is due to the fact $\lambda(g_1) \lsup{g_1}{\lambda(\overline{g_1})} = 1$ since $\lambda$ is a co-cycle.

Define the following map:
\begin{align*}
\langle \cdot, \cup \cdot\rangle: H^1(G, \hat{A}) \otimes H^3(G, A) \overset{\cup}{\longrightarrow} H^4(G, \hat{A} \otimes A) \overset{\text{(eval)}_*}{\longrightarrow} H^4(G, U(1)),
\end{align*}
where $\cup$ is the cup product and  $\text{eval}: \hat{A} \otimes A \longrightarrow U(1)$ is the evaluation map which commutes with the $G$-action. (We assume $G$ acts on $U(1)$ trivially.) To be more precise, the formula for $\langle\lambda, \cup \beta \rangle$ is given by,
\begin{align*}
\langle\lambda, \cup \beta \rangle (g_1, g_2, g_3, g_4) = \text{eval}(\lambda(g_1), \lsup{g_1}{(\beta(g_2,g_3,g_4))}).
\end{align*}

Then Equation \ref{equ:cup} means,
\begin{align}
\label{equ:cup_condition}
\langle\lambda, \cup \beta \rangle = e \in H^4(G, U(1)).
\end{align}

Thus we obtained the generalized categorical group $\CatC(G, A, \lambda, a(\omega, \beta))$, where $\omega \in H^3(G, U(1)), \beta \in H^3(G, A)$ satisfy Equation \ref{equ:cup_condition} and $a(\omega, \beta)$ is defined by Equation \ref{equ:a_omega_beta}. We will use this category in Section \ref{subsec:invariants_categorical}.

\subsection{Idempotent Completion}
\label{subsec:idempotent}
Let $\CatC$ be a category. The idempotent completion, also called Karoubi envelop or Cauchy completion, $\IC{\CatC}$ of $\CatC$ is a category defined as follows. The objects of $\IC{\CatC}$ consist of pairs $(X, \phi)$, where $X$ is an object of $\CatC$ and $\phi: X \rightarrow X$ is an idempotent, i.e., $\phi^2 = \phi$. Given two objects $(X, \phi), (X', \phi')$ of $\IC{\CatC}$,
\begin{equation*}
\Hom_{\IC{\CatC}}((X, \phi), (X', \phi')):= \phi' \circ \Hom_{\CatC}(X, X')\circ \phi = \{\psi \in \Hom_{\CatC}(X, X'):\ \phi' \psi = \psi = \psi \phi\}.
\end{equation*}
The composition of morphisms in $\IC{\CatC}$ is the same as that in $\CatC$. 

Now let $\CatC$ be the generalized categorical group $\CatC(G, A, \lambda, a)$ defined in Section \ref{subsec:generalized}. Recall that $\{P_{\chi}: \chi \in \hat{A}\}$ forms a set of complete orthogonal idempotents. It follows that the idempotents in $\Complex[A]$ are of the form
\begin{equation}
\sum\limits_{\chi \in \tilde{A}} c_{\chi} P_{\chi}, \quad c_{\chi} \in \{0,1\}.
\end{equation}
Thus there are in total $2^{|A|}$ idempotents in $\Complex[A]$. It also follows that,
\begin{align*}
\Hom((g, P_{\chi}), (g', P_{\chi'})) = \delta_{g,g'}\delta_{\chi, \chi'}\,\Complex \{P_{\chi}\},
\end{align*}
and that,
\begin{align*}
(g, P_{\chi_1} + P_{\chi_2}) \simeq (g, P_{\chi_1}) \oplus (g, P_{\chi_2}),\, \chi_1 \neq \chi_2.
\end{align*}
Therefore $\IC{\CatC}$ is a semisimple category whose non-zero simple objects are $\{(g, P_{\chi})\ : g \in G, \chi \in \tilde{A}\}$. Since the zero morphism is an idempotent, $(g, 0)$ is the zero object for any $g$. We abbreviate $(g, P_{\chi})$ as $(g, \chi)$ when no confusion arises.

Now we study the monoidal structure on $\IC{\CatC}$. Recall that for $P_{\chi_i} \in \Hom_{\CatC}(g_i,g_i), i =1,2$, we have $P_{\chi_1} \otimes P_{\chi_2} = P_{\chi_1}P_{\phi(g_1, \chi_2)} = \delta_{\phi(\overline{g_1}, \chi_1),\chi_2}P_{\chi_1}$. Then for two simple objects $(g_1, P_{\chi_1}), (g_2, P_{\chi_2})$ of $\IC{\CatC}$, define
\begin{align*}
(g_1, P_{\chi_1}) \otimes (g_2, P_{\chi_2}):= (g_1 \otimes g_2, P_{\chi_1} \otimes P_{\chi_2}) = (g_1g_2,\delta_{\phi(\overline{g_1}, \chi_1),\chi_2}P_{\chi_1}).
\end{align*}
Thus for the tensor product to be a non-zero object, $\chi_2$ must equal $\phi(\overline{g_1}, \chi_1)$. For three objects $(g_i, P_{\chi_i}), i=1,2,3$ with $\chi_2 = \phi(\overline{g_1}, \chi_1), \chi_3 = \phi(\overline{g_2}, \chi_2) =\phi(\overline{g_1g_2}, \chi_1)$, then
\begin{align*}
\{(g_1, P_{\chi_1}) \otimes (g_2, P_{\chi_2})\} \otimes (g_3, P_{\chi_3}) =(g_1g_2g_3,P_{\chi_1})=(g_1, P_{\chi_1}) \otimes \{(g_2, P_{\chi_2}) \otimes (g_3, P_{\chi_3}) \},
\end{align*}
and define the association isomorphism by
{\footnotesize 
\begin{align*}
a(g_1, g_2, g_3)_{\chi_1}P_{\chi_1}: \{(g_1, P_{\chi_1}) \otimes (g_2, P_{\chi_2})\} \otimes (g_3, P_{\chi_3}) \overset{\simeq}{\longrightarrow} (g_1, P_{\chi_1}) \otimes \{(g_2, P_{\chi_2}) \otimes (g_3, P_{\chi_3}) \},
\end{align*}}
namely, the association isomorphism is the $P_{\chi_1}$-component of $a(g_1,g_2,g_3)$ in the $\{P_{\chi}: \chi \in \hat{A}\}$ basis. If either $\chi_2$ or $\chi_3$ is not given as above, then the tensor product of the $(g_i,P_{\chi_i})\,'$s is the zero object and we define the corresponding association isomorphism as the unique zero morphism (also the identity morphism). It is direct to check the association isomorphism satisfies the Pentagon Equation.

The unit object is defined to be $(e,e) = (e, \sum_{\chi} P_{\chi}) = \oplus_{\chi} (e, P_{\chi})$. Note that $(e, P_{\chi'}) \otimes (g, P_{\chi}) = \delta_{\chi', \chi}(g, P_{\chi}) = (g, P_{\chi}) \otimes (e,P_{\phi(\bar{g}, \chi')} )$, hence the category $\IC{\CatC}$ is a multi-fusion category indexed by $\hat{A}$. Specifically, let $\IC{\CatC}_{\chi_1, \chi_2}$ be spanned additively by
\begin{align*}
\{(g, P_{\chi}): (e, P_{\chi_1}) \otimes (g, P_{\chi}) = (g, P_{\chi}) = (g, P_{\chi}) \otimes (e, P_{\chi_2}) \}.
\end{align*}
Then we have
\begin{align*}
\IC{\CatC} = \bigoplus\limits_{\chi_1, \chi_2 \in \hat{A}} \IC{\CatC}_{\chi_1,\chi_2}.
\end{align*}
and $(g, P_{\chi}) \in \IC{\CatC}_{\chi, \phi(\overline{g}, \chi)}$. For each $\chi$, $(e, P_{\chi})$ is the unit in the fusion category $\IC{\CatC}_{\chi,\chi}$. For a simple object $(g, P_{\chi})$, define the dual $(g, P_{\chi})^*:= (\overline{g}, P_{\phi(\overline{g}, \chi))}$.

We sum up the properties of $\IC{\CatC}$ as a proposition.
\begin{proposition}
Let $\CatC = \CatC(G,A, \lambda, a)$ be a generalized categorical group, then $\IC{\CatC}$ is a $\SMFC$ indexed by $\hat{A}$ where,
\begin{enumerate}
\item the simple objects correspond to elements of $G \times \tilde{A}$;
\item $(g, \chi)$ is in the sector $(\chi, \phi(\overline{g}, \chi))$;
\item $(g_1, \chi_1) \otimes (g_2, \chi_2) = \delta_{\phi(\overline{g_1}, \chi_1), \chi_2} (g_1g_2, \chi_1)$;
\item the quantum dimension of each simple object is $1$, and the dimension of each row is $|G|$, thus $\CatC$ is $\Special$.
\end{enumerate}

\end{proposition}

\subsection{Invariants from Generalized Categorical Groups}
\label{subsec:invariants_categorical}
Throughout this section, let $\CatC = \CatC(G, A, \lambda, a)$ be a generalized categorical group. We study the invariant of $3$-manifolds $Z_{\IC{\CatC}}(\cdot)$ where $\IC{\CatC}$ is the $\SMFC$ as constructed in Section \ref{subsec:idempotent}.

Let $M$ be a closed oriented $3$-manifolds and $\CatT$ be an ordered triangulation of $M$. We have $I = \hat{A}$ and $\Label{\IC{\CatC}} = G \times \hat{A}$. Recall from Section \ref{sec:main} that a $\IC{\CatC}$-coloring is a pair of maps $F = (f^0,f^1),\, f^0: \CatT^0 \longrightarrow \hat{A}, \,f^1 : \CatT^1 \longrightarrow G \times \hat{A}$ such that for any $1$-simplex $(01)$,
\begin{align*}
f^1_{01} \in \IC{\CatC}_{f^0_0,f^0_1}.
\end{align*}
Let $f^1_{01} = (g_{01}, \chi_{01})$, then
\begin{align*}
\chi_{01} = f^0_0, \quad \text{and} \quad f^0_1 = \phi(\overline{g_{01}}, f^0_0).
\end{align*}
Thus $\chi_{01}$ is uniquely determined by the coloring $f^0$. For any $2$-simplex $(012)$,
\begin{align*}
V^{\pm}_{F}(012) =
\begin{cases}
\Complex, & g_{02} = g_{01}g_{12} \\
0       , & \text{otherwise.}
\end{cases}
\end{align*}

Combing the observations, we have the following definition.
\begin{definition}
\label{def:admissible_color}
Let $\CatC, M, \CatT$ be as above. An admissible $\IC{\CatC}$-coloring is a pair of maps $\tilde{F} = (f^0,g)$, $f^0: \CatT^0 \longrightarrow \hat{A},\, g: \CatT^1\longrightarrow G$ such that,
\begin{itemize}
\item for any $1$-simplex $(01)$, $f^0_1 = \phi(\overline{g_{01}}, f^0_0)$;
\item for any $2$-simplex $(012)$, $g_{02} = g_{01}g_{12}$.
\end{itemize}
\end{definition}

Given an admissible $\IC{\CatC}$-coloring $\tilde{F} = (f^0,g)$, choose any path $p_i$ in $\CatT^1$ consisting of the edges $(v_0v_1)-(v_1v_2)- \cdots - (v_{m-1}v_m)$ where $v_0 = 0, v_m = i$, and let $g_{p_i} = g_{v_0v_1} \cdots g_{v_{m-1}v_m}$, where $g_{v_{k-1}v_k}$ is defined to be $\overline{g_{v_kv_{k-1}}}$ if $v_{k-1} > v_{k}$. It is not hard to see that $f^0_i = \phi(\overline{g_{p_i}},f^0_0)$, and thus the choice of a path $p_i$ connecting vertex $0$ to $i$ is irrelevant.

\begin{proposition}
The invariant $Z_{\IC{\CatC}}(M)$ is given by the formula,
\begin{align}
Z_{\IC{\CatC}}(M) &= \frac{1}{|G|^{|\CatT^0|}}\sum\limits_{\tilde{F}} \prod\limits_{\gamma = (ijkl) \in \CatT^3} \{a(g_{ij}, g_{jk}, g_{kl})_{f^0_{i}}\}^{\epsilon(\gamma)} \nonumber \\
                  &= \frac{1}{|G|^{|\CatT^0|}}\sum\limits_{\tilde{F}} \prod\limits_{\gamma = (ijkl) \in \CatT^3} \{a(g_{ij}, g_{jk}, g_{kl})_{\phi(\overline{g_{p_i}},f^0_0)}\}^{\epsilon(\gamma)}
\end{align}
\begin{proof}
In $\IC{\CatC}$, for any admissible $(a,b,c)$ of simple objects, one can choose $B_{ab}^c$ and $B_{c}^{ab}$ (see Section \ref{sec:main}) to be the identity map, and $\theta(a,b,c) = 1$. The quantum dimension of each simple object is $1$ and the dimension of each row is $K = |G|$.

Admissible colorings correspond to those colorings whose contribution to the summation term in Equation \ref{equ:partition} or \ref{equ:partition2} is not zero. Thus we only need to consider admissible colorings. Given an admissible coloring $\tilde{F} = (f^0,g)$, for a $3$-simplex $\gamma = (0123)$, the evaluation of Figure \ref{fig:6j} (Left) is seen to be $a(g_{01},g_{12},g_{23})_{f^0_{0}}$.

\end{proof}
\end{proposition}

\begin{corollary}
Let $\CatC = \CatC(G, A, \lambda, (\beta, \omega))$ be as defined in Section \ref{subsec:generalized}, then
\begin{align}
Z_{\IC{\CatC}}(M) &= \frac{1}{|G|^{|\CatT^0|}}\sum\limits_{\tilde{F}} \prod\limits_{\gamma = (ijkl) \in \CatT^3} \{\omega(g_{ij}, g_{jk}, g_{kl})\,f^0_{i}(\beta(g_{ij}, g_{jk}, g_{kl}))\}^{\epsilon(\gamma)}.
\end{align}
\end{corollary}

The partition function in the above corollary matches exactly the $(2+1)$-$\TQFT$ (the dual model) constructed from higher gauge theory in \cite{kapustin2013higher}, where a finite gauge group is replaced by a finite $2$-group. Thus here we provided a categorical construction of such $\TQFT$s. According to \cite{kapustin2013higher}, the $\TQFT$s thus obtained are more general than Dijkgraaf-Witten theory and provide new symmetry protected phases of matter.

\section{$2D$ Symmetry Enriched Topological Phases}
\label{sec:SET}
Symmetry plays an important role in understanding topological phases of matter. A useful approach to study topological phases is to construct exactly solvable lattice models. When anyon excitations also possess global symmetries, such a topological phase is called a symmetry enriched topological ($\SET$) phase. $\SET$s in two spacial dimension are of great interest in condensed matter physics. In \cite{cheng2016exactly}\cite{barkeshli2016reflection}\cite{heinrich2016symmetry}\cite{chang2015enriching}, exactly solvable models for a wide class of $(2D)$ bosonic $\SET$s are constructed. When the global symmetry is onsite and unitary, then the input to their models is a unitary $G$-graded fusion category, where $G$ is the global symmetry group. In this section, we show that their construction of $\SET$s extends to the framework of multi-fusion categories.

Let $\CatD = \bigoplus\limits_{g \in G}\CatD_g$ be a $G$-graded unitary fusion category and let $\tilde{\CatD}$ be the multi-fusion category obtained from $\CatD$ as given in Example \ref{example}. That is, $\tilde{\CatD} = \bigoplus\limits_{g,h \in G}\tilde{\CatD}_{g,h}$ is indexed by $G$ where $\tilde{\CatD}_{g,h} = \CatD_{\bar{g}h}$. The tensor products in $\tilde{\CatD}$ are the same as those in $\CatD$ and for any $g \in G$ the unit $\unit_g$ in $\tilde{\CatD}_{g,g}$ is the unit $\unit$ in $\CatD_{e}$ (and also the unit in $\CatD$). $\tilde{\CatD}$ is spherical since $\CatD$ is spherical. (Unitarity implies sphericity.) Also, for any $g \in G$, $K(\tilde{\CatD}_{g}) = \sum\limits_{h} K(\tilde{\CatD}_{g,h}) = \sum\limits_{h} K(\CatD_{g^{-1}h}) = K(\CatD)$. Thus, $\tilde{\CatD}$ is a $\Special$ $\SMFC$.



Assume $\CatD$ is multiplicity free. As in Section \ref{sec:main}, for any admissible $(a,b,c)$ of simple objects, we choose a basis element $B_{c}^{ab} \in  \Hom(c, a\otimes b)$ and $B_{ab}^c \in \Hom(a \otimes b, c)$ such that,
\begin{align*}
\langle B_{c}^{ab} , B_{ab}^c \rangle = \Tr(B_{ab}^cB_{c}^{ab}) = \theta(a,b,c),
\end{align*}
where $\theta(a,b,c) = \sqrt{d_ad_bd_c}$.

Let $M$ be an oriented $3$-manifold and $\CatT$ be a triangulation of $M$. If $M$ has no boundary, then the partition function $Z_{\tilde{\CatD}}(M)$ is given by Equation \ref{equ:partition} or Equation \ref{equ:partition2} as a state-sum model. By definition, a $\tilde{\CatD}$-coloring $F = (f^0,f^1)$ assigns to each vertex ordered by $k$ a group element $f^0_k \in G$ and assigns to each $1$-simplex $(ij)$ a simple object $f^1_{ij} \in \tilde{\CatD}_{f^0_if^0_j} = \CatD_{\overline{f^0_i}f^0_j}$. It is direct to check the partition function $Z_{\tilde{\CatD}}(M)$ thus obtained is the same as the one given in \cite{barkeshli2016reflection}. More generally, when $M$ is bounded by a surface $\partial M$, the wave function associated with $M$ is defined by
\begin{align}
\Psi(\partial M, F) =  \sum\limits_{\tilde{F}: \tilde{F}_{| \partial M} = F} \prod\limits_{\tau \in \CatT^0} K^{-1} \prod\limits_{\alpha \in \CatT^1} d_{f^1_{\alpha}} \,  \prod\limits_{\beta \in \CatT^2} \theta_F(\beta)^{-1} \, \prod\limits_{\gamma \in \CatT^3} \tilde{Z}^{\epsilon(\gamma)}_F(\gamma;B),
\end{align}
where $F$ is a coloring of $\CatT$ restricted to $\partial M$  and the summation on the right hand side is over all colorings $\tilde{F}$ extending $F$.

For $g \in G$, let $\lsup{g}{F}$ be the coloring $(g. f^0, f^1)$, namely, the color on each vertex is multiplied on the left by $g$ while the color on each edge remains unaltered. Clearly, $\lsup{g}{F}$ is a well-defined color and that $\Psi(\partial M, F) = \Psi(\partial M, \lsup{g}{F})$. 

\vspace{1cm}
\noindent\textbf{Acknowledgment} $\quad$ The first author acknowledges the support from the Simons Foundation and would like to thank Meng Cheng and Ryan Thorngren for helpful discussions. The second author is partially supported by NSF grant DMS-1411212.

\bibliographystyle{plain}
\bibliography{MultiFusionbib}

\end{document}